\documentclass[12pt,draftcls,onecolumn]{IEEEtran} 
\usepackage{cite}
\usepackage{graphicx,color}
\usepackage{multirow}
\usepackage{float}
\usepackage[cmex10]{amsmath}
\usepackage{amssymb}
\usepackage{subfigure}
\usepackage{bm}
\usepackage{color}
\usepackage{relsize}

\usepackage[top=2.5cm, bottom=3.0cm, left=2.5cm, right=2.5cm]{geometry}

\everymath{\displaystyle}

\newcommand{\beqnum}{\begin{equation}\begin{array}{lcl}}
\newcommand{\eeqnum}{\end{array}\end{equation}}
\newcommand{\beqnom}{\begin{eqnarray}}
\newcommand{\eeqnom}{\end{eqnarray}}
\newcommand{\beqnc}{\begin{center}\begin{eqnarray}}
\newcommand{\eeqnc}{\end{eqnarray}\end{center}}
\newcommand{\beqnlm}{\begin{equation}\vspace{-.5cm}\begin{array}{lll}}
\newcommand{\eeqnlm}{\end{array}\end{equation}}\vspace{-.5cm}
\newcommand{\beq}{\begin{eqnarray*}}
\newcommand{\eeq}{\end{eqnarray*}}
\newcommand{\bef}{\begin{figure}[tbh!]}
\newcommand{\enf}{\end{figure}}

\newtheorem{define}{\bf Definition}

\newcommand{\vep}{\varepsilon}
\newcommand{\lb}{\lambda}

\newcommand{\R}{\mathbb{R}}
\newcommand{\lf}{\left\lfloor}
\newcommand{\rr}{\right\rceil}

\newtheorem{montheo}{\bf Theorem}
\newtheorem{rem}{\bf Remark}

\newtheorem{proposition}{\bf Proposition}

\usepackage{graphics} % for pdf, bitmapped graphics files
\usepackage{epsfig} % for postscript graphics files

\title{\LARGE \bf
Higher order super-twisting algorithm for perturbed chains of integrators
of arbitrary order
\thanks{This research was partially supported by the iCODE Institute, research project of the IDEX Paris-Saclay, and by the Hadamard Mathematics LabEx (LMH) through the grant number ANR-11-LABX-0056-LMH in the ``Programme des Investissements d'Avenir''.}
}

\author{Yacine Chitour, Mohamed Harmouche, Salah Laghrouche% <-this % stops a space
\thanks{Y. C. is with L2S, Universite Paris 11, CNRS 91192 Gif-sur-Yvette, France. {\tt\small yacine.chitour@lss.supelec.fr}}%
\thanks{M. H. is with Actility, Paris, France.
        {\tt\small mohamed.harmouche@actility.com}}%
\thanks{S. L. is with OPERA Laboratory, UTBM, Belfort, France.
{\tt\small salah.laghrouche@utbm.fr}}%
}

\begin{document}

\maketitle
\pagenumbering{arabic}

%%%%%%%%%%%%%%%%%%%%%%%%%%%%%%%%%%%%%%%%%%%%%%%%%%%%%%%%%%%%%%%%%%%%%%%%%%%%%%%%
\begin{abstract}
\noindent
In this paper, we present a generalization of the super-twisting algorithm for perturbed chains of integrators of arbitrary order. This Higher Order Super-Twisting (HOST) controller, which extends the approach of Moreno and als., is homegeneous with respect to a family of dilations and can be continuous. Its design is derived from a first result obtained for pure chains of integrators, the latter relying on a geometric condition introduced by the authors. The complete result is established using a homogeneous strict Lyapunov function which is explicitely constructed.
%Our result is compared with the approach presented with \cite{Kamal2014}, where the authors present a quasicontinuous form of Higr Order SuperTwisting without proof. 
The effectiveness of the controller is finally illustrated with simulations for a chain of integrator of order four, first pure then perturbed, where we compare the performances of two HOST controllers.\end{abstract}

%\tableofcontents
%%%%%%%%%%%%%%%%%%%%%%%%%%%%%%%%%%%%%%%%%%%%%%%%%%%%%%%%%%%%%%%%%%%%%%%%%%%%%%%%
\section{Introduction}

The control of nonlinear industrial systems is a challenging task because these systems are difficult to characterize and suffer from parametric uncertainty. Parametric uncertainty arises from various operating conditions and external perturbations that affect the physical characteristics of systems. This must be considered during the control design so that the controller counteracts the effect of variations and guarantees performance under different operating conditions. Sliding Mode Control (SMC) \cite{Utkin,Slotine1984,Levant93} is well-known for control of nonlinear systems and renowned for its insensitivity to bounded parametric uncertainty and external disturbance. This technique is based on applying discontinuous control on a system which ensures convergence of the output function (sliding variable) in finite time to a manifold of the state-space, called the sliding manifold. In practice, SMC suffers from \emph{chattering}; the phenomenon of finite-frequency, finite-amplitude oscillations in the output which appear because the high-frequency switching excites unmodeled dynamics of the closed loop system \cite{Utkin1999}.  Higher Order Sliding Mode Control (HOSMC) is an effective method for chattering attenuation \cite{Emelyanove}, where the discontinuous control is applied on a higher time derivative of the sliding variable. In this way, the sliding variable and also its higher time derivatives converge to the origin. As the discontinuous control does not act upon the system input directly, chattering is automatically reduced.

Many HOSMC algorithms exist in contemporary literature for control of nonlinear systems with bounded uncertainty. These algorithms are robust because they preserve the insensitivity of classical sliding mode, and maintain the performance characteristics of the closed loop system. Levant for example, has presented a method of designing arbitrary order sliding mode controllers for Single Input Single Output (SISO) systems in \cite{Levant2001}. Laghrouche et al. \cite{Laghrouche2} have proposed a two part integral sliding mode based control to deal with the finite time stabilization problem and uncertainty rejection problem separately. Dinuzzo et al. have proposed another method in \cite{Dinuzzo09}, where the problem of HOSMC has been treated as Robust Fuller's problem. Defoort et al. \cite{Defoort2009} have developed a robust multi-input multi-output HOSMC controller, using a constructive algorithm with geometric homogeneity based finite time stabilization of an integrator chain. Harmouche et al. have presented their homogeneous controller in \cite{HARMOUCHE1} based on the work of Hong \cite{Hong}. Sliding mode with homogeneity approach was also used in \cite{Levant2003,Levant2005}, to demonstrate finite time stabilization of the arbitrary order sliding mode controllers for SISO systems \cite{Levant2001}. A Lyapunov-based approach for arbitrary HOSMC controller design was presented in \cite{Harmouche_CDC12, Laghrouche_CST}. In this work, it was shown that a class of homogeneous controllers that satisfies certain conditions, could be used to stabilize perturbed integrator chains.  

The main drawback of these controllers is that they produce a discontinuous control signal \cite{Kamal2014},  at least at the origin. In order to build continuous controllers and still have finite time convergence, a standard trick consists in using a relative degree extension as advocated for instance in \cite{Levant_Springer}, namely to consider the extra equation $\dot u=v$, where $u$  is a control and $v$ is a virtual HOSM discontinuous control input. This controller should use  the output and its the first $r$ derivatives, where $r$ is the sliding mode order. However, to complete that procedure, it is necessary to know the bounds of the uncertainties and their first derivatives and, more restrictively, to suppose that the time integral of $v$ is uniformly bounded for all initial conditions (see Remark 2). 
To overcome this problem and to also get finite time convergence, Kamal et al. \cite{Kamal2014} propose a generalization of the well-known continuous super-twisting algorithm  for high order relative degree system with respect to the output (cf.\cite{Levant97}),  ensuring finite time convergence of the sliding variable and its $r$ first derivatives to zero, by using a continuous control signal and only  information about the sliding variables and its $r-1$ derivatives \cite{Kamal2016}, \cite{Fridman_2016}. The convergence conditions and Lyapunov analysis have been only given up to order three and a higher order controller is just suggested without proof. Other results solving this problem were proposed in \cite{Basin2015} and \cite{Edwards2016}. These algorithms are not homogeneous and thus they can not achieve the $r$-th order of sliding precision with respect to the sliding variable \cite{Levant2003}. 

In the present paper, we provide a homogenous  HOST controller for arbitrary order with a complete argument as well as other HOST controllers. %which are just continuous at the origin. 
Our analysis is based on the use of a homogeneous strict Lyapunov function %and homogeneity arguments 
for an extended system. The resulting HOST controllers are either continuous or just continuous at the origin and both ensure finite-time stabilization, first for a pure chain of integrators and then for a perturbed one. 

To describe our results, recall that a perturbed chain of integrators of length $r$ reads $\dot z_i=z_{i+1}$ for $1\leq i\leq r-1$ and $\dot z_r=\gamma u+\varphi$ where $\gamma$ is a positive measurable signal lower and upper bounded with known positive constants and both 
$\dot \gamma$ and $\dot \varphi$ are bounded by known positive constants. Note that we do not assume that the additive perturbation $\varphi$ is bounded. We first provide a HOST controller for a pure chain of integrators (i.e.,  $\gamma$ is constant and $\varphi$ is equal to zero) based on standard controllers for a pure chain of integrators verifying in addition a geometric condition. The convergence proof using the HOST controller relies on the existence of a homogeneous strict Lyapunov function $W$ associated with an extended system. The construction of $W$ is explicit once the standard controllers for the pure chain of integrators are given. We then prove these HOST controllers can be used for perturbed chains of integrators. In the case $\gamma$ constant (let say equal to one), one must recall that one can stabilize in finite time such a perturbed chain of integrators without using HOST controllers. Indeed, by setting $z_{r+1}:=u+\varphi$ and $v:=\dot u$, one gets $\dot z_{r+1}=v+\dot \varphi$ and a perturbed chain of length $r+1$ is stabilized with the control $v$ and the bounded uncertainty $\dot\varphi$. This can be done at the price of increasing the length of the chain of integrators (i.e., the relative degree of the output), which can be a serious drawback in some applications. This technics is referred to as \emph{extension of relative degree}, cf.\cite{Levant_Springer, Laghrouche2}. One should therefore see our HOST solution as an alternative to the increase of the length of the chain of integrators. In the case of non constant $\gamma$ and assuming that $0\leq \gamma_m\leq \vert \gamma\vert \leq \gamma_M$, we provide a solution for arbitrary length under a smallness condition on $1-\gamma_m/\gamma_M$. The controller we propose depends on both the state and time and presents a finite number of discontinuities in the time variable only and becomes continuous if  $\varphi$ is assumed to be bounded. Note that 
we are able to tune the parameters of our controller in the general case of non constant $\gamma$ in terms of the homogeneous Lyapunov function $W$.

The paper is organized as follows. 
%The first part of Section~\ref{section2} gathers the main definitions regarding differential inclusions, homogeneity and finite time convergence for homogeneous differential inclusions. In the second part of the section, we treat in detail the one-dimensional case, reinterpreting the solution proposed by \cite{Moreno} while describing our strategy for higher order. 
In Section~\ref{section3}, we present our results for the general case, first addressing the stabilization by HOST of a pure chain of integrators of arbitrary order and then explaining how to generalize to a perturbed chain by means of homogeneity arguments. We prove two results, one in the case of constant $\gamma$ and the second one for non constant $\gamma$. We close Section~\ref{section3} by providing explicit examples of standard controllers for pure chains of integrators which do verify the required geometric condition. 
We finally demonstrate in Section~\ref{section4} the efficiency of our HOST algorithm for a perturbed chain of integrators of order four.

{\bf Notations.}
In this paper, we use $\R$ and $\vert\cdot\vert$ to denote the set of real numbers and a fixed norm on $\R^r$ respectively, where $r$ is a positive integer. For $\lb>0$, let $D_\lb$ be the $r\times r$  matrix defined by $\hbox{diag}(\lb^r,\cdots,\lb)$.
For $m$ positive integer, let $e_1,\cdots,e_m$ and $J_m$ denote the canonical basis of $\R^m$ and the $m$-th Jordan block as $(J_m)_{ij}=\delta_{i,j+1}$, $1\leq i,j\leq m$, respectively. 
If $M$ is a subset of $\R^r$, we use $\overline{M}$ to denote its closure. If $x\in\mathbb{R}$, we denote by $[x]$ the integer part of $x$ i.e., the smallest integer not greater than $x$. 
We define the function $sign$ as the multivalued function defined on $\mathbb{R}$ by $sign(x)=x/\vert x\vert$ for $x\neq 0$ and $sign(0)=[-1,1]$. Similarly, for every $a\geq 0$ and $x\in \mathbb{R}$, we use $\lf x\rr^a$ to denote $\left| x \right|^a sign(x)$. Note that $\lf \cdot\rr^a$ is a  continuous function for $a>0$ and is of class $C^1$ with derivative equal to $a\left| \cdot \right|^{a-1}$ for $a\geq 1$. If $V:\R^r\rightarrow \R^p$ is a differentiable mapping, we use $\partial_jV$ to denote the partial derivative of $V$ with respect to the $j$-th coordinate $z_j$ 
and more generally $\partial_\xi V$ if $V$ depends on a scalar coordinate $\xi$.
We refer to \cite{Levant_Springer, Bernuau2014} for the definitions of {\it a Filippov
differential inclusion}, {\it Homogeneity} and  {\it Asymptotic and Finite time stability}. %In particular, we will use repeatedly \cite{Levant_Springer, Bernuau2014} where it is shown that, for a differential inclusion $\dot z\in F(z)$, $z\in\R^r$ with negative degree, asymptotic and finite time stability are equivalent. 

\section{Higher order super-twisting continuous feedback for a chain of integrator}\label{section3}
%In this section, we essentially generalize what has been done for the first-order integrator. 
The strategy consists first in building an appropriate feedback 
for a pure chain of integrator and then in tackling the perturbed case by a homogeneity argument.
\subsection{Stabilization of a pure chain of integrator of arbitrarily order}\label{sec:pure}
%\beqnum \label{pure_int}
%	\left\{
%		\begin{array}{ccl}
%			\dot z_1 &=& z_2,\\
%				&\vdots&\\
%			\dot z_r &=& u,
%		\end{array}
%	\right.
%\eeqnum
\begin{define}\label{def0} {\it Let $r$ be a positive integer. The $r$-th order chain of integrator
$(CI)_r$ is the single-input control system given by 
$$
(CI)_r\ \ \ \dot z=J_rz+u e_r,\ \ \ z=(z_1,\cdots,z_r)^T\in \R^r,\ \ u\in \R.
$$
For $\kappa<0$ and $p>0$ with $p+(r+1)\kappa\geq0$, set $p_i := p + (i-1)\kappa, \ 1\leq i\leq r+1$. For $\varepsilon>0$, let $\delta_\vep:\R^r\rightarrow \R^r$ and $\psi_\vep:\R^{r+1}\rightarrow \R^{r+1}$ be the family of dilations associated with $\left( p_1, \cdots , p_r \right)$ and $\left( p_1, \cdots , p_{r+1} \right)$ respectively.}
\end{define}
In the spirit of \cite{Laghrouche_CST, Harmouche_CDC12}, we put forwards geometric conditions on certain stabilizing feedbacks $u_0(\cdot)$ for $(CI)_r$ and corresponding Lyapunov functions $V_1$. These conditions will be instrumental for building a super-twisting feedback law as shown in the following theorem.
%{\color{green}
%!!!!!!!!!!!!!!!!!!!!!!!!!!!!!!!!!!!!!!
%
%Cette reference concerne le CDC harmouche, baghdadi, laghrouche ainsi que nos papier avec la condition geometrique
%
%!!!!!!!!!!!!!!!!!!!!!!!!!!!!!!!
%}
\begin{montheo}\label{theo1}
Let $r$ be a positive integer. Assume that there exists a feedback law 
$u_0:\R^r\rightarrow\R$, homogeneous  with respect to $(\delta_\vep)_{\vep>0}$ of degree $p_{r+1}$ such that the closed-loop system $(CI)_r$ with $u_0$ is finite time stable and the following conditions hold true:
\begin{description}
	\item[$(i)$] the function $z\mapsto J_rz+u_0(z)e_r$ is homogeneous of degree $\kappa$ with respect to $(\delta_\vep)_{\vep>0}$ and there exists a continuous positive definite function $V_1:\R^r\rightarrow \R_+$, $C^1$ except at the origin, homogeneous with respect to $(\delta_\vep)_{\vep>0}$ of degree $2p_{r+1}$ such that there exists $c>0$ and $\alpha\in (0,1)$ for which the time derivative of $V_1$ along non trivial trajectories of $(CI)_r$ verifies $\dot V_1 \le -cV_1^\alpha$.
		\item[$(ii)$] $z\mapsto \partial_r V_1(z)$ is homogeneous of non negative degree with respect to $(\delta_\vep)_{\vep>0}$ and
$z\mapsto \partial_r V_1(z)u_0(z)$ is non negative over $\mathbb{R}^r$.
\end{description}
Then, for every $k_P\geq 1$ and $k_I>0$,  $(CI)_r$ is stabilized in finite time by the HOST controller
\beqnum\label{ST-feedback}
		u_{ST}(z,t) = k_P u_0(z) -k_I \int_{0}^{t} \partial_r V_1(z(s))ds.
\eeqnum
Moreover, write the (time-varying) closed loop system $\dot z=J_rz+u_{ST}e_r$ as the differential inclusion over $\R^{r+1}$
\beqnum\label{pure-1}
 \dot z=J_rz+(k_Pu_0+\xi)e_r,\ \  	\dot \xi = -k_I \partial_r V_1,
\eeqnum
where $\xi:=-k_I \int_{0}^{t}\partial_r V_1dt$. 
Then, there exists $A>0$ so that the continuous function $W:\R^{r+1}\rightarrow \R$, defined by
\beqnum\label{ST-LF}
	W(z,\xi) = A\big(V_1(z)+ \xi^2/2 k_I\big)^{2-\alpha}-z_r\xi,
\eeqnum
is positive definite, $C^1$ except at the origin, homogeneous with respect to $(\psi_\vep)_{\vep>0}$ and there exists $d>0$ for which the time derivative of $W$ along non trivial trajectories of \eqref{pure-1} verifies $\dot W \le -dW^{1/(2-\alpha)}$. As a consequence, trajectories of \eqref{pure-1} converge to zero in finite time.
\end{montheo}

{\bf Proof of Theorem~\ref{theo1}.}
%The closed loop system $\dot z=J_rz+u_{ST}e_r$ can be written as
%\beqnum\label{pure-1}
% $$
% \dot z=J_rz+(k_Pu_0+\xi)e_r,\ \  	\dot \xi = -k_I \partial_r V_1,
%%	
%%	
%%	\left\{
%%		\begin{array}{ccl}
%%			\dot z_1 &=& z_2,\\
%%				&\vdots&\\
%%			\dot z_r &=&  k_P u_0   +  \xi,\\
%%			\dot \xi &=& -k_I \partial_r V_1,
%%		\end{array}
%%	\right.
%\eeqnum
%where
%$\xi:=-k_I \int_{0}^{t}\partial_r V_1dt$. 
First of all, notice that $\alpha$ must be equal to $1+\kappa/(2p_{r+1})\geq 1/2$ and 
$\partial_rV_1$ is homogeneous with respect to $(\delta_\vep)_{\vep>0}$ of degree $p+(r+1)\kappa\geq 0$.
 
Consider first the positive definite function $V:=V_1 +  \xi^2/2 k_I$, which is homogeneous with respect to $(\psi_\vep)_{\vep>0}$ of degree $2p_{r+1}$. 
Using Items $(i)$ and  $(ii)$, the time derivative of $V$ along trajectories of \eqref{pure-1} verifies
%\begin{align*}
%\dot V&=\sum_{i=1}^{r-1}\partial_i V_1z_{i+1}+\partial_r V_1( k_P u_0   +  \xi)
%-\xi \partial_r V_1
%=\dot V_1+(k_P-1)\partial_r V_1u_0\leq -cV_1^{\alpha}(z).
%\end{align*}
$\dot V=\dot V_1+(k_P-1)\partial_r V_1u_0\leq -cV_1^{\alpha}$.

For $A>0$, the function $W$ defined in \eqref{ST-LF} is continuous, $C^1$ except at the origin, homogeneous with respect to $(\psi_\vep)_{\vep>0}$ of degree $2(2-\alpha)p_{r+1}=p_r+p_{r+1}$. Since $V$ is positive definite, $W$ is also positive definite for $A$ large enough.  
The time derivative of $W$ along trajectories of \eqref{pure-1} verifies
\begin{align*}
\dot W&\leq -c(2-\alpha) AV_1-k_Pu_0\xi+k_Iz_r\partial_r V_1-\xi^2.
\end{align*}
For $A$ large enough, the right-hand side of the previous inequality is negative definite and smaller than $-dW^{1/(2-\alpha)}$ for some positive constant $d$. Therefore, trajectories of \eqref{pure-1} converge to zero in finite time and hence
$(CI)_r$ is stabilized in finite time by the feedback law $u_{ST}$. $\blacksquare$
\begin{rem}\label{rem:deg1}{\it If $p+(r+1)\kappa=0$, then $\partial_rV_1$ is of zero homogeneity degree with respect to $(\delta_\vep)_{\vep>0}$ and $\alpha=1/2$. In that case, $W$
is homogeneous with respect to $(\psi_\vep)_{\vep>0}$ of degree $3p_{r+1}$.
Notice that the constants $A$ and $d$ can be explicitely computed if $u_0$ and $V_1$ are explicitely given.}
%Moreover, by taking $\overline{W}:=W^{5/6}$, one gets a continuous positive-definite function  $C^1$ except at the origin, homogeneous with respect to $(\psi_\vep)_{\vep>0}$ and its time derivative along non trivial trajectories of \eqref{pure-1} verifies
%$\dot {\overline{W}}\leq -e \overline{W}^{1/2}$ for some positive constant $e$.}
\end{rem}
% \hspace*{\fill} $\blacksquare$ 
%\begin{rem}In the above argument, we only need %$u_0$ to be continuous at zero %(implying that $u_0(0)=0$) 
%$\partial_r V_1$ to be bounded in an open neighborhood of the origin in order to get the following proposition.\end{rem}
%We also have the following proposition.
%\begin{proposition}\label{Prop1}
%Assume that the hypotheses of Theorem~\ref{theo1} are satisfied with $p+(r+1)\kappa=0$.
%%and $\partial_r V_1$ is bounded in a neighborhood of the origin. 
%Then, there exists $c>0$ and a continuous positive-definite function $W:\R^{r+1}\rightarrow \R_+$, $C^1$ except at the origin, homogeneous with respect to $(\psi_\vep)_{\vep>0}$ and the time derivative of $W$ along non trivial trajectories of \eqref{pure-1} verifies
%%\beqnum\label{1DD-der}
%$\dot W\leq -c W^{1/2}$.
%%\eeqnum
%%for some positive constant $c$. 
%\end{proposition}
%{\bf Proof of Proposition~\ref{Prop1}.} One checks that \cite[Theorem 4.1]{Bernuau2014} applies. $\blacksquare$ 
%Indeed, one clearly verifies that the right-hand side of \eqref{pure-1} is a locally essentially bounded vector field giving rise to a multivalued function $F$ which is non-empty, compact and convex. Therefore $F$ satisfies the standard assumptions through the Filippov regularization procedure (cf. \cite{Filippov} for precise results) and is in addition homogeneous of degree $\kappa$ with respect tothe family of dilations $\psi_\vep$ defined on $\R^{r+1}$ by 
%\hspace*{\fill} $\blacksquare$ 
%}
%====================================================================
\subsection{Stabilization of a perturbed chain of integrators: case of $\gamma$ constant }
In this subsection, we apply the previous results to get finite-time convergence of the perturbed chain of integrators defined next by 
\beqnum \label{perturbed_int}
\dot z=J_rz+( \gamma u + \varphi)e_r
%	\left\{
%		\begin{array}{ccl}
%			\dot z_1 &=& z_2,\\
%				&\vdots&\\
%			\dot z_r &=& \gamma u + \varphi, 
%		\end{array}
%	\right.
\eeqnum
where the time-varying functions $\gamma(\cdot)$ and $\varphi(\cdot)$ are globally Lipschitz  over $\R_+$ and
verify the following: there exist $\gamma_m,\gamma_M>0$ and 
$\overline{\gamma},\overline{\varphi}\geq 0$ such that, for every $t\geq 0$ it holds
\beqnum \label{bound-varphi}
0<\gamma_m\leq \gamma(t)\leq\gamma_M,\quad \vert\dot \gamma(t)\vert\leq \overline{\gamma},\quad
\vert\dot \varphi(t)\vert\leq \overline{\varphi}.
\eeqnum
%\begin{align}
%0<\gamma_m&\leq \gamma(t)\leq\gamma_M, \label{bound-gamma}\\
%\vert\dot \gamma(t)\vert\leq \overline{\gamma}&, \ 
%\vert\dot \varphi(t)\vert\leq \overline{\varphi}.\label{bound-varphi}
%\end{align}
%In particular, notice that both $\gamma$ and $\varphi$ are globally Lipschitz with Lipschitz constants upper bounded by $\overline{\gamma}$ and $\overline{\varphi}$ respectively. 

We now want to derive conditions under which the super-twisting feedback defined in Eq.~\eqref{ST-feedback} stabilizes \eqref{perturbed_int} in finite time.
%As mentioned in Introduction, we restrict in this paper to the case where $\gamma\equiv 1$ and 
We obtain the following theorem. 

\begin{montheo}\label{theo2}
Consider the perturbed chain  of integrators defined by \eqref{perturbed_int}, where the time-varying function $\gamma(\cdot)$ and $\varphi(\cdot)$ verify $\gamma\equiv \gamma_m$
%\eqref{bound-gamma} 
and \eqref{bound-varphi} respectively. Assume that there exists a continuous homogeneous feedback law 
$u_0$ and a Lyapunov function $V_1$ verifying the assumptions $(i)$ and $(ii)$
of Theorem~\ref{theo1} with $p+(r+1)\kappa=0$. Then, for every positive gains $k_P\geq 1$ and $k_I$, there exists $\lambda>0$ only depending on the gains and $\overline{\varphi}$ such that 
the feedback law $u_{ST}(D_\lb z,\lb t)/\gamma_m$, where %the matrix $D_\lb$ is defined in Eq. \eqref{matriceDl} and 
$u_{ST}$ is given in \eqref{ST-feedback}, stabilizes \eqref{perturbed_int} in finite-time.
In particular, $u_{ST}(D_\lb z,\lb t)/\gamma_m$ is continuous.
\end{montheo}
{\bf Proof of Theorem~\ref{theo2}.}
Fix now some $k_P\geq 1$ and $k_I>0$. Associate to every absolutely continuous function $z:\mathbb{R}_+\rightarrow\mathbb{R}^r$ the following function for $t\geq 0$
$\xi(t):=-k_I\int_0^t\partial_r V_1(z(s))ds+\varphi(t)$.
%g(t):=-\frac{\gamma(t)-\gamma_0}{\gamma_0}u_{ST}(z(t),t).
%\big(k_Pu_0(z)-k_I\int_0^t\frac{\partial V_1}{\partial z_r}(z(s))ds\big).
The closed-loop system  \eqref{perturbed_int}  with $u_{ST}$ can be written as 
\beqnum \label{perturbed_int_1}
\dot z=J_rz+(k_pu_0(z)+\xi)e_r,\ \dot \xi=-k_I\partial_r V_1+ \dot \varphi(t).
%	\left\{
%		\begin{array}{ccl}
%			\dot z_1 &=& z_2,\\
%				&\vdots&\\
%			\dot z_r &=&  k_Pu_0(z)   + \xi,\\
%			\dot \xi &=& -k_I\frac{\partial V_1}{\partial z_r}+ \dot \varphi(t).
%		\end{array}
%	\right.
\eeqnum
%Note that one has, for every trajectory $z(\cdot)$ of \eqref{perturbed_int_1} and every $t\geq 0$, 
%\beq\label{est-gamma}
%\vert\frac{\gamma(t)-\gamma_0}{\gamma_0}\vert\leq \frac{\overline{\gamma}-\underline{\gamma}}{\overline{\gamma}+\underline{\gamma}}:=\rho<1.
%\eeq
Clearly, \eqref{perturbed_int_1} corresponds to the differential inclusion~\eqref{pure-1} perturbed by  the time-varying vector field over $\R^{r+1}$ given by $(0,\cdots,0,\dot\varphi(t))^T$ or, equivalently, by the multifunction $(0,\cdots,0,[-\overline{\varphi},\overline{\varphi}])^T$ taking values in the subsets of $\R^{r+1}$.
Let $W$ be the Lyapunov function defined in \eqref{ST-LF}. Along non trivial trajectories of \eqref{perturbed_int_1}, one gets, for every $t\geq 0$, that
\begin{align}\label{estW}
\dot W\leq -dW^{2/3}+\partial_\xi W\dot\varphi(t)\leq -dW^{2/3}+\overline{\varphi}\vert \partial_\xi W\vert.
\end{align}
%associated with \eqref{pure-1} furnished by Remark~\ref{rem:deg1}. Along non trivial trajectories of \eqref{perturbed_int_1}, one gets, for every $t\geq 0$, that
%\begin{align}\label{estW}
%\dot W\leq -cW^{1/2}+\partial_\xi W\dot\varphi(t)\leq -cW^{1/2}+\overline{\varphi}\vert \partial_\xi W\vert.
%\end{align}
According to Remark~\ref{rem:deg1}, the homogeneity degree of $\vert \partial_\xi W\vert$ is equal to $3p_{r+1}-p_{r+1}=2p_{r+1}$, i.e., the homogeneity degree of $W^{2/3}$. 
%Indeed, if $F:\R^{r+1}\rightrightarrows\R^{r+1}$ denotes the differential inclusion defined by \eqref{pure-1}, then $\max_{y\in F(z,\xi)}\nabla W(z,\xi)y\leq -dW^{2/3}$
%for $(z,\xi)\in \R^{r+1}$ by Remark~\ref{rem:deg1}. Moreover, 
%%by eventually reducing $1/2$ in the previous inequality
%by using the argument in \cite[Theorem 7.1]{Bernstein2005} adapted for homogeneous differential inclusions, there exists $e>0$ such that, for every $(z,\xi)\in \R^{r+1}$, one has 
%$-eW^{2/3}\leq \max_{y\in F(z,\xi)}\nabla W(z,\xi)y$.
%We now rewrite $\nabla W(z,\xi)y$ for $y\in F(z,\xi)$ as
%$$
%\nabla W(z,\xi)y=\sum_{i=1}^{r-1}\partial_i Wz_{i+1}+\partial_r W
%(k_Pu_0+\xi)-k_I\partial_\xi W\partial_r V_1.
%$$
%By a standard homogeneity argument, %every term in the previous sum is a homogeneous function with the same degree with respect tothe family of dilations $(\delta_\vep(z),\vep^{p_{r+1}}\xi)$ and thus 
%$\Big\vert \partial_\xi W\Big\vert$ has the same homogeneous degree as $W^{1/2}$ with respect to$(\psi_\vep)_{\vep>0}$. 
One deduces from Eq.~\eqref{estW} that there exists $\varphi_*>0$ such that 
$\dot W\leq -\frac{d}2W^{1/2}$,
along trajectories of System \eqref{perturbed_int_1} if $\overline{\varphi}\leq \varphi_*$. We thus have proved the theorem under the previous restriction on $\overline{\varphi}$. 
To remove it, consider the standard time-coordinate change of variable along trajectories of \eqref{perturbed_int} defined, for every $\lb>0$, by
%\beqnum\label{TC}
$y(t)=D_\lb z(t/\lb)$.
%\eeqnum
Under the hypotheses of the theorem, \eqref{perturbed_int} can be rewritten
$\dot y=J_ry+(u_\lb/\gamma_m+\varphi_\lb)e_r$, 
%\beqnum \label{perturbed_int_y}
%	\left\{
%		\begin{array}{ccl}
%			\dot y_1 &=& y_2,\\
%				&\vdots&\\
%			\dot y_r &=& \frac{u_\lb}{\gamma_m}+\varphi_\lb,%\frac{u(\frac {t}{\lb}) + \varphi(\frac {t}{\lb})}{\lb}.
%		\end{array}
%	\right.
%\eeqnum
where one has set, for $t\geq 0$, $u_\lb(t):=u(t/\lb)$ and $\varphi_\lb(t):=\varphi(t/\lb)$. 
Note that, for a. e. $t\geq 0$, 
$\vert \dot\varphi_\lb\vert \leq \overline{\varphi}/\lb$.
By taking $\lb\geq \overline{\varphi}/\varphi_*$, one gets that $
\vert\dot\varphi_\lb\vert \leq \varphi_*$. We apply the previous stabilization result and conclude.
$\blacksquare$ 

Notice that the choice of $\varphi_*$ can be made explicit once $u_0$ and $V_1$ are explicitely given.
%Let's recall the following Lemma:
%\begin{lemme}[Corollary 1 of \cite{Levant2005}] The global uniform finite-time stability of homogeneous differential equations (Filippov inclusions) with
%negative homogeneous degree is robust with respect tohomogeneous perturbations causing locally small changes of the equation (inclusion) graph.
%\end{lemme} 
%
%Based on the previous Lemma, we get automatically the following result
%\begin{lemme}
%There exist a small positive constant $\bar \phi \le \left| \partial_r V_1  \right|$, such that for all $\left | \dot \varphi \right| \le \bar \phi$, the following system perturbed system remains finite time stable with the feedback $u$.
%\end{lemme}
%Clearly the previous Lemma proves the existence of parameter $\bar \phi$. Otherwise we can get an estimation through simulation.\\
%In our future research we try to estimate  $\bar \phi$ based on strong Lyapunov function.
%\begin{rem}
%	In the previous Lemma, we get the utility of the condition $[iiii]$.
%\end{rem}
%%%%%%%%%%%%%%%%%%%%%%%%%%%%%%%%%%%%%%%%%%%%%%%%%%%%%%%%%%%%%%%%%%%%%%%%%%%%%%%%%%%%%%%%%%%%%%%%%%%%%%%%%%%
%{\color{blue}
\subsection{Stabilization of a perturbed chain of integrators: general case}
In this subsection, we apply Subsection~\ref{sec:pure} to get finite-time convergence 
of \eqref{perturbed_int}
where the time-varying functions $\gamma(\cdot)$ and $\varphi(\cdot)$ are measurable over $\R_+$ and
verify \eqref{bound-varphi}.

\begin{rem}\label{rem:erd}
Note that the technics of extension of relative degree seems difficult to be implemented in case $\gamma(\cdot)$ is not constant. Indeed, 
if we set again $z _{r+1}:= \varphi + \gamma u$, one gets $\dot z_{r+1} = \dot \varphi + \dot\gamma u + \gamma \dot u$, which can be written as $\dot z_{r+1} = \gamma v + \Phi$, with $\dot u = v$ and $\Phi = \dot \varphi + \dot\gamma u$. In order to pursue with the classical arguments, it is necessary to have $\Phi$ bounded, which occurs only if $u$ itself is.  One is therefore lead to determine a feedback law $v$ for a perturbed $r$-th chain of integrators such that the time integral of $v$ is uniformly bounded for all initial conditions. 
\end{rem}

We now want to derive conditions under which the super-twisting feedback defined in Eq.~\eqref{ST-feedback} stabilizes System \eqref{perturbed_int} in finite time.
%As mentioned in Introduction, we restrict in this paper to the case where $\gamma\equiv 1$ and w
We obtain the following theorem. 

\begin{montheo}\label{theo3}
Consider the perturbed chain  of integrators defined by \eqref{perturbed_int}, where the time-varying function $\gamma(\cdot)$ and $\varphi(\cdot)$ verify \eqref{bound-varphi}.
%\eqref{bound-gamma} 
Assume that there exists a continuous feedback law 
$u_0$ and a Lyapunov function $V_1$ verifying the assumptions 
of Theorem~\ref{theo1} with $p+(r+1)\kappa=0$.
Set $\gamma_d=(\gamma_M+\gamma_m)/2$ and $\delta_\gamma=1-\gamma_m/\gamma_M\in (0,1)$.
Then, for every gains $k_P\geq 1$ and $k_I>0$, there exists $\delta_0\in (0,1)$ and $\lambda_0>0$ 
only depending on the gains and the constants in  Eq~\eqref{bound-varphi} such that, if $\delta_\gamma\leq \delta_0$
and $\lambda\geq \lambda_0$, \eqref{perturbed_int} is stabilized in finite-time by
the feedback law 
$$
u_{st}(z,t)=\frac1{\gamma_d}(k_P u_0(D_{\lb} z)-k_I%\mathlarger{\int_{ [\lambda t]}^{\lambda t}}\partial_r V_1(D_\lb z(s))ds)$.
\int_{ [\lambda t]}^{\lambda t}\partial_r V_1(D_\lb z(s))ds).
$$
%where the matrix $D_\lb$ is defined in Eq. \eqref{matriceDl} stabilizes \eqref{perturbed_int} in finite-time. 
\end{montheo}

{\bf Proof of Theorem~\ref{theo3}.} As in the proof of Theorem~\ref{theo2}, fix $k_P\geq 1$ and $k_I>0$. For $t\geq 0$, define 
%\begin{align}\label{eq:discont}
$\xi(t):=-\gamma(t)\gamma_d/k_I%\mathlarger{\int_{[t]}^t}\partial_r V_1(z(s))ds+\varphi(t)$.
\int_{[t]}^t\partial_r V_1(z(s))ds+\varphi(t)$, which yields an absolutely continuous function on each interval $[n,n+1]$, with a discontinuity at $t=n$. %g(t):=-\frac{\gamma(t)-\gamma_0}{\gamma_0}u_{ST}(z(t),t).
%\big(k_Pu_0(z)-k_I\int_0^t\frac{\partial V_1}{\partial z_r}(z(s))ds\big).
%*\end{align}
%for some positve constant $C$ to be chosen later.
System \eqref{perturbed_int} with $u_{st}$ in \eqref{perturbed_int}  can be written, a.e. on each interval $[n,n+1]$, 
as 
\beqnum \label{perturbed_int_3}
	\left\{
		\begin{array}{ccl}
			\dot z&=&J_rz+(k_Pu_0(z)   + \xi+(\gamma(t)-\gamma_d)k_Pu_0/\gamma_d)e_r,\\
			\dot \xi &=& -k_I\partial_r V_1(z)+ \dot \varphi(t)
			-k_I\frac{\gamma(t)-\gamma_d}{\gamma_d}\partial_r V_1(z)
			-k_I\frac{\dot \gamma(t)}{\gamma_d}\mathlarger{\int_{[t]}^t}\partial_r V_1(z(s))ds.
		\end{array}
	\right.
\eeqnum
To pursue the argument, set $v_1:=\max_{z\in \mathbb{R}^r\setminus\{0\}}\vert\partial_r V_1(z)\vert$ and notice that
$$
\forall t\geq 0,\ \ \vert (\gamma(t)-\gamma_d/\gamma_d\vert\leq \delta_\gamma/(2-\delta_\gamma)
\leq \delta_\gamma,\quad
\vert  \int_{[t]}^t\partial_r V_1(z(s))ds\vert\leq 
\max_{s\in [[t],t]}\vert \partial_r V_1(z(s))\vert.
%v_1&=\max_{z\in \mathbb{R}^r\setminus\{0\}}\left\vert \partial_r V_1(z)\right\vert.
%\end{align*}
$$
Consider now the Lyapunov function $W$ defined in \eqref{ST-LF}. Along non trivial trajectories of System~\eqref{perturbed_int_3}, one gets 
a.e. on each interval $[n,n+1]$ with $n$ integer, that 
\begin{align}
\dot W&\leq -dW^{2/3}+\delta_\gamma k_P\vert \partial_r Wu_0(z)\vert+
\vert \partial_\xi W\vert(\overline{\varphi}+ k_Iv_1(\delta_\gamma +
\overline{\gamma}/\gamma_d)).\label{est-W-2}
\end{align}
Exactly as for $\vert \partial_\xi W\vert$, $\vert \partial_r Wu_0(z)\vert$ is homogeneous with respect to $(\psi_\vep)_{\vep>0}$ of degree $2p_{r+1}$. One deduces that there exist $\delta_0$, $\varphi_0>0$ and $\tilde{\gamma}$ such that, if 
%\beqnum\label{restrictions}
$\delta_\gamma \leq \delta_0$, $\overline{\varphi}\leq \varphi_*$, $\overline{\gamma}/\gamma_d\leq \tilde{\gamma}$,
%\eeqnum 
then one has, for every $t\geq 0$, 
%\beqnum\label{est-W-1}
$\dot W\leq -\frac{d}2W^{2/3}$,
%\eeqnum
along non trivial trajectories of \eqref{perturbed_int_3}.
Assuming the above restrictions on the bounds, we prove the theorem. One cannot conclude immediately as in the proof of Theorem~\ref{theo2} since $\xi(\cdot)$ is not continuous and thus $t\mapsto W(z(t),\xi(t))$ is as well discontinuous at integer valued times along trajectories of 
\eqref{perturbed_int_3}. To address that issue, we introduce some notations: for $n$ positive integer, set 
$
Y_n=\lim_{t\rightarrow n,\ t<n}W(z(t),\xi(t))\quad Z_n=\lim_{t\rightarrow n,\ t>n}W(z(t),\xi(t)).
$
Dividing $\dot W\leq -\frac{d}2W^{2/3}$ by $3W^{2/3}(t)$ and integrating, one gets, for every positive integer $n$ and $t\in (n,n+1)$ so that the trajectory of \eqref{perturbed_int_3} remains non trivial, that
$Z_n>W(t)>Y_{n+1}$ and $Y_{n+1}^{1/3}-Z_n^{1/3}\leq -d/6$.
%\vert Z_n-Y_n\vert&=\left\vert  \mathlarger{\int_{n}^{n+1}}\frac{\partial V_1}{\partial z_r}(z(s))ds \right\vert\leq v_1,\\
%-dW^{1/2}(t)&\leq \dot W(t)\leq -\frac{c}2W^{1/2}(t).
Moreover, to estimate the jump of $W$ at discontinuity times, we notice that
$
Z_{n+1}-Y_{n+1}=W(z(n+1),\xi(n+1))-W(z(n+1),\xi_{n+1}^-),
$
where 
\begin{align*}
\vert \xi_{n+1}^--\xi(n+1)\vert&=\vert\gamma(n+1)k_I/\gamma_d\int_n^{n+1}\partial_r V_1(z(s))ds\vert\leq D_0\tilde{\gamma},
\end{align*}
with $D_0>0$ only depending on the gains $k_P$ and $k_I$.
Because $\vert \partial_\xi W\vert$ is homogeneous with respect to $(\psi_\vep)_{\vep>0}$ of degree $2p_{r+1}$, we deduce that there exists $D_1>0$ only depending on the gains $k_P$ and $k_I$ such that, for every non negative integer $n$ so that the trajectory of \eqref{perturbed_int_3} remains non trivial on $(n,n+1)$, one gets
$$
\vert Z_{n+1}-Y_{n+1}\vert=\vert \int_{\xi(n+1)}^{\xi_n^-}\partial_\xi W(z(n+1),\eta)d\eta\vert\leq D_1\tilde{\gamma}\max(Z_{n+1}^{2/3},Y_{n+1}^{2/3}).
$$
By dividing the previous inequality by $Z_{n+1}^{2/3}+Z_{n+1}^{1/3}Y_{n+1}^{1/3}+Y_{n+1}^{2/3}$, we deduce that
$\vert Z_{n+1}^{1/3}-Y_{n+1}^{1/3}\vert\leq D_1\tilde{\gamma}$. 
By taking $\tilde{\gamma}$ smaller than $d/12D_1$ and using the above results, we get for $n\geq 1$ that
\begin{equation}\label{finite1}
Y_{n+1}^{1/3}-Y_n^{1/3}=(Y_{n+1}^{1/3}-Z_n^{1/3})+ (Z_n^{1/3}-Y_n^{1/3})\leq -d/12.
\end{equation}
The finite-time convergence to the origin of non trivial trajectories of \eqref{perturbed_int_3} follows at once.

To remove most of these restrictions, we proceed as in the proof of Theorem~\ref{theo2}, i.e., by considering the time-coordinate transformation $y(t)=D_\lb z(t/\lb)$. For $t\geq 0$ and $\lambda>0$, set 
$u_\lb(t)=u(t/\lb)$, $\gamma_\lb(t)=\gamma(t/\lb)$ and $\varphi_\lb(t)=\varphi(t/\lb)$.
Then, $\vert \dot\gamma_\lb\vert \leq \overline{\gamma}/\lb$ and $\vert \dot\varphi_\lb\vert \leq \overline{\varphi}/\lb$.
By taking $\lb\geq \lambda_0:=\max(\overline{\varphi}/\varphi_*,1/\gamma_d\tilde{\gamma})$, one gets that $
\vert \dot\varphi_\lb\vert \leq \varphi_*$ and $\vert\dot \gamma_\lb/\gamma_d\vert\leq \tilde{\gamma}$. We can now apply the previous stabilization result and conclude. Note though that we are not able with this trick to remove the restriction on $\delta_d$. $\blacksquare$ 

Notice that the choices of $\tilde{\gamma}$ and $\varphi_*$ can be made explicit once $u_0$ and $V_1$ are explicitely given. If we impose an extra restriction on $\varphi$, we get a continuous HOST feedback. This is explained in the following property.
\begin{proposition}\label{extra}
Consider the same hypotheses as in Theorem~\ref{theo3} and, in addition suppose that 
$\vert \varphi\vert\leq \varphi_M$ for some known non negative constant $\varphi_M$. 
Then the same conclusion as in Theorem~\ref{theo3} is reached with the feedback law 
$u_{st}$ defined by 
\beqnum\label{feedback-per}
u_{st}(z,t)=k_P u_0(D_{\lb} z)-k_I\int_{0}^{\lambda t}\partial_r V_1(D_\lb z(s))ds)/\gamma_d.
\eeqnum
In particular, the feedback $u_{st}$ is continuous.
\end{proposition}

{\bf Proof of Proposition~\ref{extra}.} 
We define now the integral variable $\xi$ as follows for $t\geq 0$
\begin{equation}\label{rem2}
\xi(t):=-k_I\int_{0}^t\partial_r V_1(z(s))ds+\gamma_d\varphi(t)/\gamma(t).
\end{equation}
%for some positve constant $C$ to be chosen later.
The closed-loop \eqref{perturbed_int}  together with $u_{st}$ can be written as 
\beqnum \label{perturbed_int_4}
	\left\{
		\begin{array}{ccl}
			\dot z&=&J_rz+(k_Pu_0(z)   + \xi+(\gamma(t)-\gamma_d)/\gamma_d(k_Pu_0(z)+\xi))e_r,\\
			\dot \xi &=& -k_I\partial_r V_1+\gamma_d
			(\dot\varphi \gamma-\varphi\dot \gamma)/\gamma^2
					\end{array}
	\right.
\eeqnum
Eq.~\eqref{est-W-2} becomes
$\dot W\leq -dW^{2/3}+\delta_\gamma \vert \partial_r W\vert \vert k_Pu_0(z)+\xi\vert+
\vert \partial_\xi W\vert \gamma_d(\overline{\varphi}\gamma_M+\varphi_M\overline{\gamma})/\gamma_m^2$.
From the above inequality, one finishes the argument as in the proof of Theorem~\ref{theo3}. $\blacksquare$ 

\begin{rem}\label{rem:bdd}
From Eqs.~\eqref{feedback-per} and \eqref{rem2}, one deduces that, for $t\geq 0$,  
$u_{st}(z,t)$ is equal to $k_P u_0(D_{\lb} z)-\gamma_d\xi(t)-\varphi(t)/\gamma(t)$
 and remains bounded along a given trajectory.
\end{rem}
%%%%%%%%%%%%%%%%%%%%%%%%%%%%%%%%%%%%%%%%%%%%%%%%%%%%%%%%%%%%%%%%%%%%%%%%%%%%%%%%%%%%%%%%%%%%%%%%%%%%%%%%%%%%%%%%%%%%%%%%%%%%%%%%%%%%%%%%
\subsection{Feedbacks $u_0$ and Lyapunov functions $V_1$ verifying the assumptions of Theorems~\ref{theo2} and ~\ref{theo3}}
We next provide examples of controllers $u_0$ and Lyapunov functions $V_1$ 
satisfying the conditions of Theorem~\ref{theo2} for $r\geq 2$. We next assume that $p=1$ and thus $\kappa=-1/(r+1)$.
\subsubsection{Hong's controller}\label{Hong0}
Such a controller is simply borrowed from \cite{Hong}. In that reference, 
the convergence is proved by using a Lyapunov function $V_0$ explicitly constructed for that purpose. The latter function does not
match the the assumptions of Theorem~\ref{theo2} and we have to modify it to get the required Lyapunov function $V_1$.
%The controller provided in \cite{Hong} is defined as follows. 
Let $l_1,\cdots,l_r$ positive real numbers. We define, for $ i=0,...,r+1$, the functions 
$v_0\equiv 0$ and for $1\leq i\leq r-1$, $v_{i+1} = -l_{i+1} \lfloor\lfloor z_{i+1} \rceil^{\beta_i } -\lfloor v_i \rceil^{\beta_i } \rceil^{\alpha_{i+1}/\beta_i}$,
%\beqnum\label{Hongfunction}
%\left\{\begin{array}{ccl}
%	v_0&=&0,\\
%	v_{i+1} &=& -l_{i+1} \lfloor\lfloor z_{i+1} \rceil^{\beta_i } -\lfloor v_i \rceil^{\beta_i } \rceil^{\alpha_{i+1}/\beta_i},\\
%							&& \quad i=1,\cdots r-2\\
%	v_r  &=& -l_r  \left\lfloor \lfloor z_{r} \rceil^{\beta_{r-1} } -\lfloor v_i \rceil^{\beta_{r-1} } \right \rceil^{\alpha_{r}/\beta_{r-1}},
%\end{array}\right.
%\eeqnum
where $\beta_0 =p_2$ and $(\beta_i + 1)p_{i+1} = \beta_0 + 1$, and $\alpha_i=p_{i+1}/p_i$, for $i=1,...,r$.
%, and
%\beqnum\label{beta00}
%	\beta_0 =p_2, \ (\beta_i + 1)p_{i+1} = \beta_0 + 1 > 0, \ i=1,...,r-1.\\
%\eeqnum
One then takes the controller $u_0$ to be equal to the continuous function $v_r$. Consider the Lyapunov function $V_0= \sum_{i=1}^{r} W_i$, where the positive real-valued functions $W_i$, $1\leq i\leq r$ are given by $W_i = \int_{v_{i-1}}^{z_i}{w_i(z_1,\cdots,z_{i-1},s) ds}$
with $w_i=\lfloor z_{i} \rceil^{\beta_{i-1} } -\lfloor v_{i-1} \rceil^{\beta_{i-1 } } , \quad i = 1,\cdots , r$.
Then, there exists $l,l_1,\cdots,l_r>0$ such that the time derivative of $V$ along every non trivial trajectory of $\dot z=J_rz+u_0e_r$ satisfies
$\dot V_0 \le-  l \bar V_0 ^{\frac{2+2\kappa}{2 + \kappa}}$. This proves that $u_0$ stabilizes $\dot z=J_rz+ue_r$ to the origin in finite-time. Let $\lambda:=2/(2r-1)<1$. Note that 
$1-\lambda=\beta_{r-1}/(1+\beta_{r-1})$.
Take now $V_1= V_0 ^\lambda/\lambda$. A simple computation yields $\partial_r V_1=(\lfloor z_{r} \rceil^{\beta_{r-1} } -\lfloor v_{r-1} \rceil^{\beta_{r-1 }})/V_0^{\frac{\beta_{r-1}}{\beta_{r-1} + 1 }}$.
%
%\beqnum
%	\frac{\partial V_1}{ \partial z_r } &=& V^{\lambda-1}\frac{\partial V}{ \partial z_r } = V^{\lambda-1}\frac{\partial W_r}{ \partial z_r }
%	=\frac{ \lfloor z_{r} \rceil^{\beta_{r-1} } -\lfloor v_{r-1} \rceil^{\beta_{i-1 }} }{V^{\frac{\beta_{r-1}}{\beta_{r-1} + 1 }}}.\\
%			%&=&\frac{ \left| \lfloor z_{i} \rceil^{\beta_{r-1} } -\lfloor v_{r-1} \rceil^{\beta_{i-1 }} \right| }{V_1^{\frac{\beta_{r-1}}{\beta_{r-1} + 1 }}} sign\left( z_r - v_{r-1} \right),\\
%\eeqnum
One then checks that  $\partial_r V_1$ is homogeneous of degree zero with respect to  $\delta_\vep$, globally bounded, and continuous except at the origin.

\subsubsection{Modified Hong's Controller}\label{Hong1}
The following controller is a hybrid form between the continuous controller presented by Hong \cite{Hong} and a terminal sliding mode approach also presented by Hong et al. in \cite{Hong2004}. Note that its form is very close to the controller proposed in \cite{Hong2005_1}. Let $\kappa$, the $l_i$'s, $\alpha_i$'s and the $\beta_i$'s as before. The functions $v_i$ are defined as above for $1\leq i\leq r-1$ but we now set $v_r=-l_rw_r$ with $w_r:= \left( \left| z_r \right| ^{\beta_{r-1}}  + \left| \bar v_{r-1} \right| ^{\beta_{r-1}  }  \right)^{\alpha_{r}/\beta_{r-1}}$ $sign\left( z_r - v_{r-1} 	\right)$.
%
%Let $k=\frac{-1}{r+1}$ and $l_1,\cdots,l_r$ positive real numbers. We define, for $ i=0,...,r-1$, $p_i=1+(i-1)k$
%\beqnum
%		\bar v_0&=&0,\\
%		\bar v_{i+1} &=& l_{i+1}  \left|  \left| z_{i+1} \right|^{\beta_i }  + \left| \bar v_i \right|^{\beta_i }  \right|^{\alpha_{i+1}/\beta_i}, i=1,\cdots r-2,\\
%\eeqnum
%where $\alpha_i=\frac{p_{i+1}}{p_i}$, for $i=1,...,r$, and the $\beta_i$'s verify Eq.~\eqref{beta00}. Then, for $ i=0,...,r-1$, we define
%\beqnum
%\left\{\begin{array}{ccl}
%	v_0&=&0,\\
%	v_{i+1} &=& -l_{i+1} \lfloor\lfloor z_{i+1} \rceil^{\beta_i } -\lfloor v_i \rceil^{\beta_i } \rceil^{\alpha_{i+1}/\beta_i}, i=1,\cdots r-2\\
%	%v_r  &=& -l_r \left( \sum_{j=1}^{r}   \left| z_j \right| ^{\beta_i }  \right)^{\alpha_{i+1}/\beta_i} sign\left( z_r - v_{r-1} \right).\\
%	v_r  &=& -l_r \left( \left| z_r \right| ^{\beta_{r-1}}  + \left| \bar v_{r-1} \right| ^{\beta_{r-1}  }  \right)^{\alpha_{r}/\beta_{r-1}} sign\left( z_r - v_{r-1} 	\right).
%\end{array}\right.
%\eeqnum
The controller $u_0$ is then taken equal to  $v_r$ and it stabilizes $\dot z=J_r z+ue_r$ in finite-time. To see that, consider the positive definite function $V_0= \sum_{i=1}^{r} W_i$ where, for $i=1,...,r$, one has $W_i = \int_{v_{i-1}}^{z_i}{w_i(z_1,\cdots,z_{i-1},s) ds}$ with
$w_i$, $1,\cdots , r-1$ defined as in Hong's controller and $w_r$ defined just above. 
%
%\beq
%	w_i &=& \lfloor z_{i} \rceil^{\beta_{i-1} } -\lfloor v_{i-1} \rceil^{\beta_{i-1 } } , \quad i = 1,\cdots , r-1,\\
%	w_r &=&  \left( \left| z_r \right| ^{\beta_{r-1}}  + \left| \bar v_{r-1} \right| ^{\beta_{r-1}  }  \right) sign\left( z_r - v_{r-1}\right).
%\eeq
%Then we get
%\beqnum
%	W_i &=&	\frac{ \left| z_i  \right|^{\beta_{i-1} + 1} - \left| v_{i-1}  \right|^{\beta_{i-1} + 1}}{\beta_{i-1} + 1} - \left\lfloor v_{i-1}\right\rceil ^{\beta_{i-1}} \left(  z_i - v_{i-1} \right), \quad i=1,\cdots r-1,\\
%	W_r &=& \frac{1}{\beta_{r-1} + 1} \left|  \left\lfloor z_r \right\rceil^{\beta_{r-1} + 1}   -   \left\lfloor v_{r-1} \right\rceil^{\beta_{r-1} + 1}    \right|
%	        + \left|\bar v_{r-1} \right| ^{\beta_{r-1}} \left| z_r  - v_{r-1}   \right|
%	\eeqnum
%
 One deduces from an argument entirely similar to that of \cite{Hong} that there exists $l>0$ such that, one has along trajectories of the closed-loop system $\dot V_0\le -  l V_0^{\frac{2+2\kappa}{2 + \kappa}}$.
Finite-time convergence to the origin follows immediately. Finally remark that the feedback control law $v_r$ is continuous at zero. One can then apply the results given in Section~\ref{section3}. The actual Lyapunov function $V_1$ is again taken of the form $V_0^\lambda/\lambda $ with $\lambda=1/(1+\beta_{r-1})$. A simple computation yields $\partial_r V_1=V_0^{\lambda-1}w_r$ which is 
homogeneous of degree zero with respect to $(\delta_\vep)_{\vep>0}$.

%%%%%%%%%%%%%%%%%%%%%%%%%%%%%%%%%%%%%%%%%%%%%%%%%%%%%%%%%%%%%%%%%%%%%%%%%%%%%%%%%%%%%%%%%%%%%%%%%%%%%%%%%%%

\section{Simulations}\label{section4}
In this section, we verify the effectiveness of our design through simulations.
We deal with a chain of integrator of order four and we show the robustness with respect to perturbations. Consider the fourth order integrator system given by 
$\dot z=J_4z+(\gamma u+\varphi(t))e_4$.
%\beqnum
%	\dot z_1 &=& z_2,\\
%	\dot z_2 &=& z_3,\\
%	\dot z_3 &=& z_4,\\
%	\dot z_4 &=& u + \varphi(t).
%\eeqnum
We study in the following subsections two cases:  
the first one deals with $\varphi\equiv 0$ and $\gamma \equiv 1$, and the second case considers $\varphi\ne 0$ and $\gamma \ne 1$, both for the continuous and discontinuous HOST controllers corresponding to the Hong's controller $u_0^H$ and to the modified Hong's controller $u_0^{MH}$ and respectively.
For all subsequent simulations, the control parameters are tuned as follows:
$l_1=l_2=1$, $l_3 = 4$, $l_4 = 8$ and $\kappa = -1/5$ with initial condition 
$z_1(0) = -5$, $z_2(0) = 2$ and $z_3(0) = z_4(0) = 4$.
%%%%%%%%%%%%%%%%%%%%%%%%%%%%%%%%%%%%%%%%%%% 
\subsection{Simulation of pure integrator chain for $u=u_0$}
\begin{figure}[htbp!]
\centering
\subfigure[control law $u$]{
    \includegraphics[width= 7 cm]{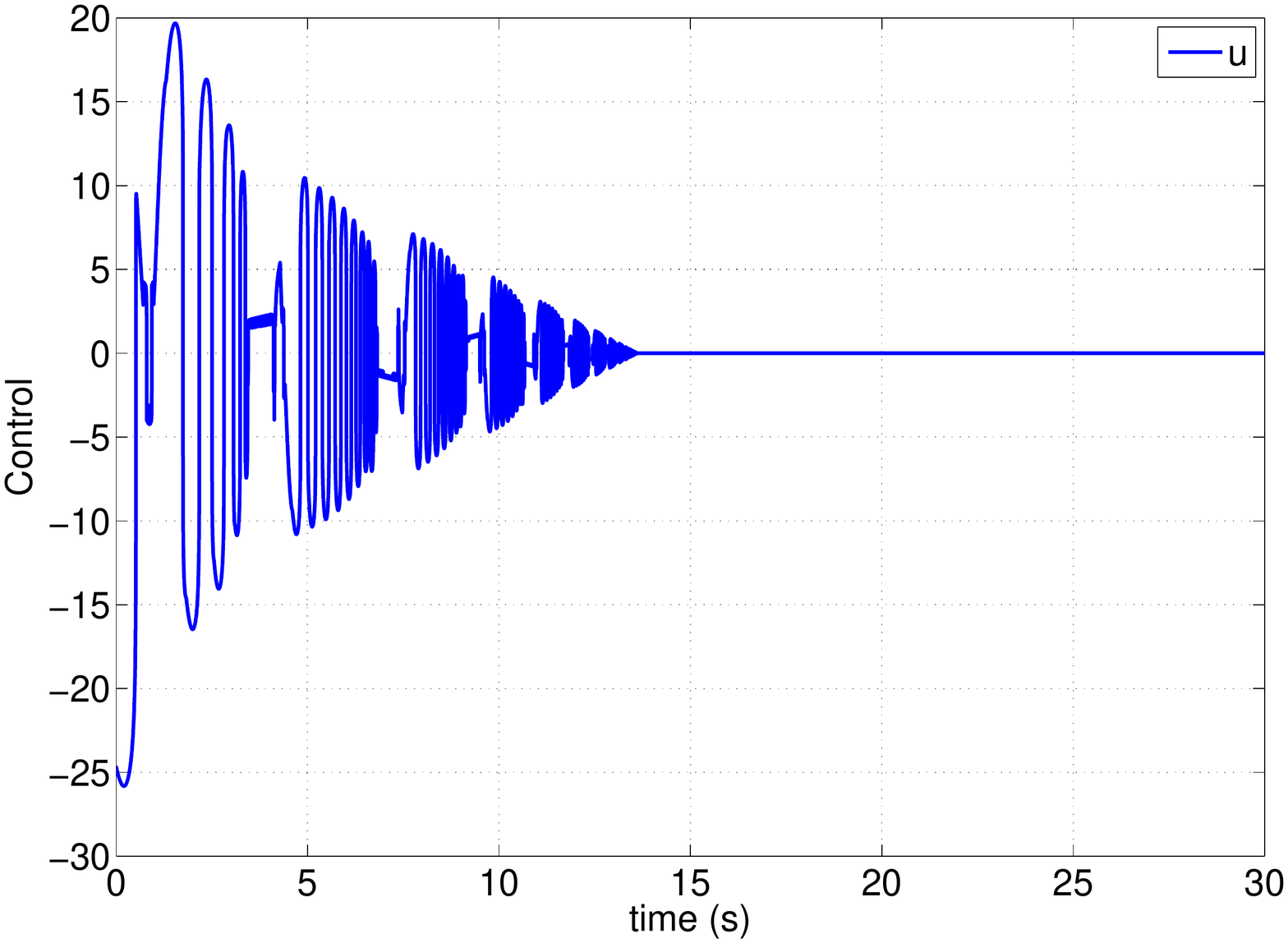}
    \label{u0}
}\hspace{-1cm}
\subfigure[$z_1$, $z_2$, $z_3$ and $z_4$]{
    \includegraphics[width= 7 cm]{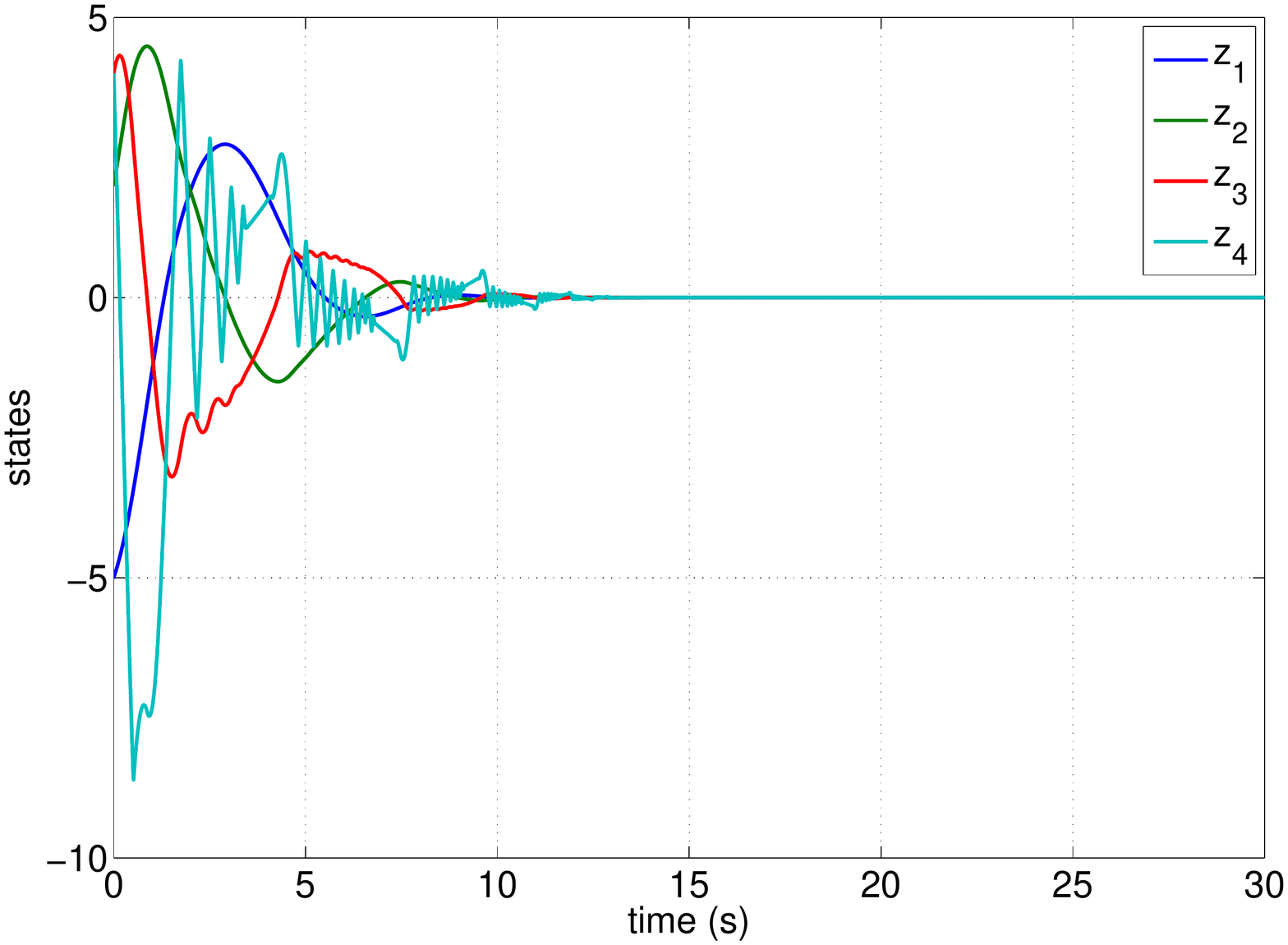}
    \label{s0}
}
\caption{Pure integrator chain without integration action (Hong's controller)}
\label{s0_u0}
\end{figure}

\begin{figure}[htbp!]
\centering
\subfigure[control law $u$ versus time ($s$).]{
    \includegraphics[width= 7 cm, height = 3.4 cm]{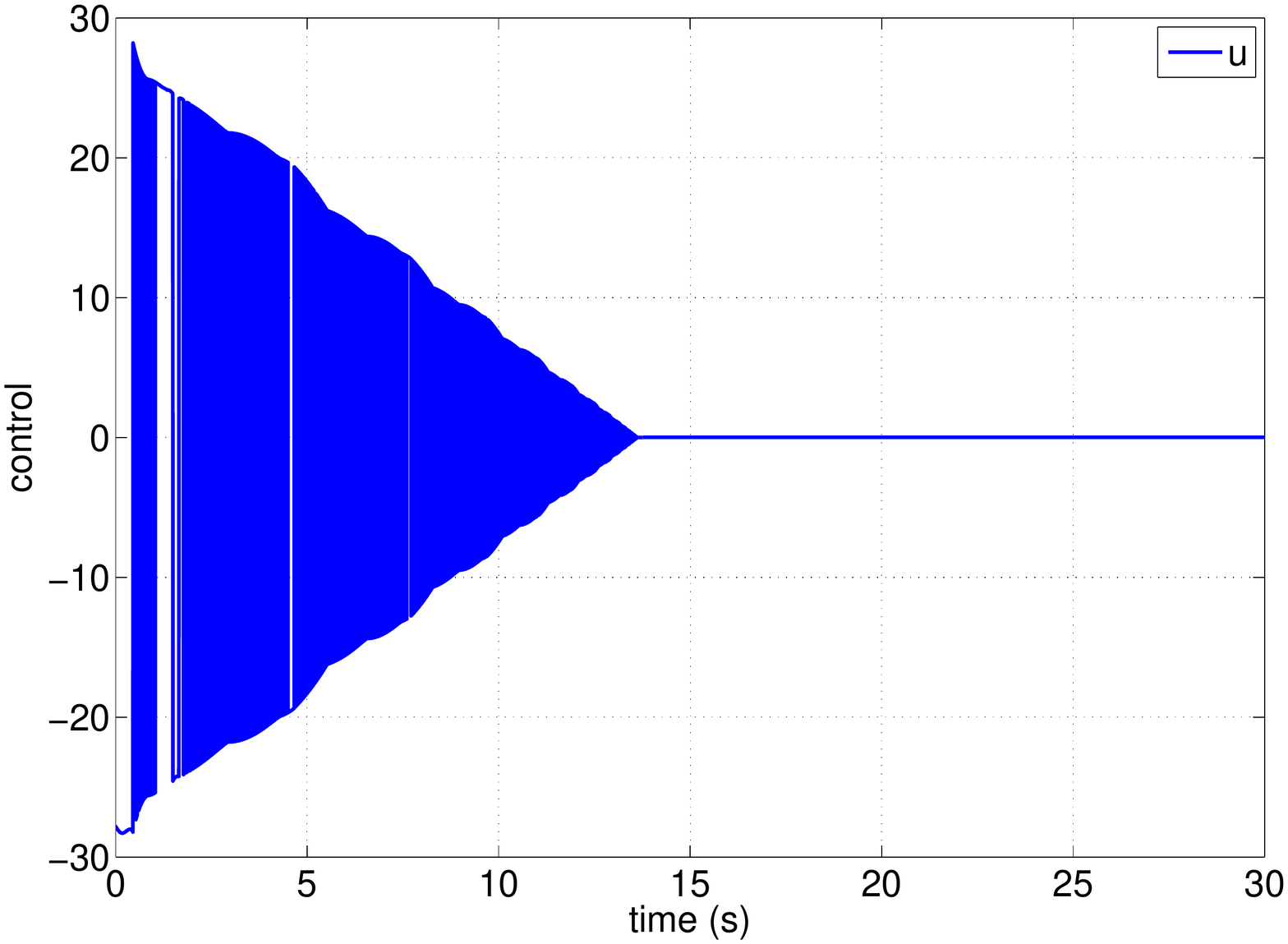}
    \label{u0_qq}
}\hspace{-1cm}
\subfigure[$z_1$, $z_2$, $z_3$ and $z_4$]{
    \includegraphics[width= 7 cm, height = 3.4 cm]{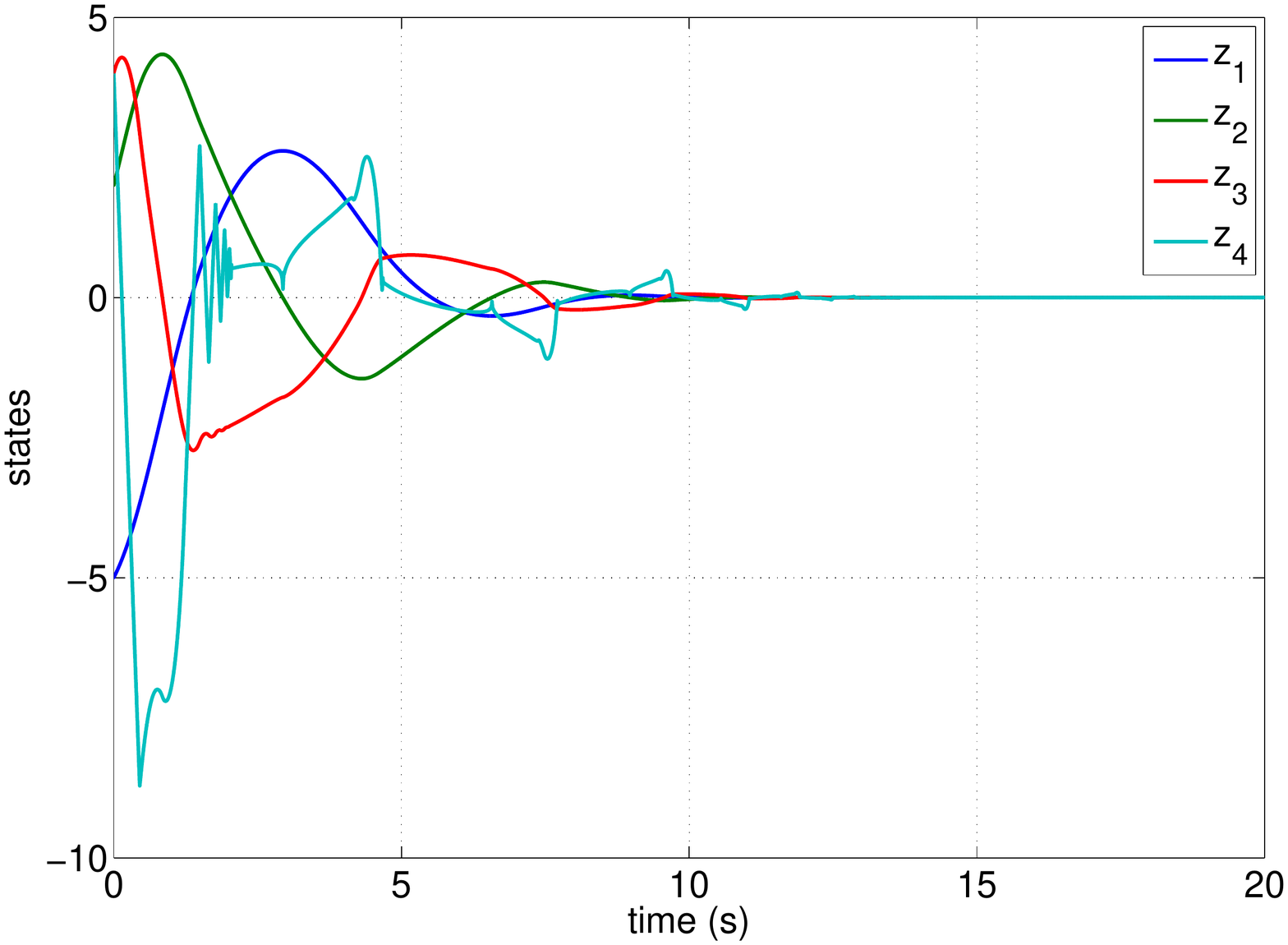}
    \label{s0_qq}
}
\caption{Pure integrator chain without integration action (Modified Hong's controller)}
\label{s0_qq_u0}
\end{figure}

We start by stabilizing the pure integrator chain, (i.e., with $\varphi\equiv 0$) by the controller $u = u_0$ where $u_0$ represents either the Hong's controller $u_0^H$ or the modified Hong's controller $u_0^{MH}$  given Sections \ref{Hong0} and \ref{Hong1} respectively. Figures \ref{u0} and \ref{u0_qq} represents $u_0^H$ and $u_0^{MH}$ presented in Sections \ref{Hong0} and \ref{Hong1} respectively. These controllers force the state $(z_1,\ z_2,\ z_3,\ z_4)$ to converge to zero in finite time, as shown in Figure \ref{s0} and Figure \ref{s0_qq_u0}.\\

%%%%%%%%

%%%%%%%%%%%%%%%%%%%%%%%%%%%%%%%%%%%%%%%%%%%%%%%%%%%%%%%%%%%%%%%%%%%%%%%%%%%%%%%%%%%%%%%%%%%%%%%%%%%%%%%%%%%%%%%%%
\subsection{Stabilisation of pure integrator chain by HOST - $\varphi\equiv 0$}
\begin{figure}[htbp!]
\centering
\subfigure[Control law $u$]{
    \includegraphics[width= 7 cm, height = 3.4 cm]{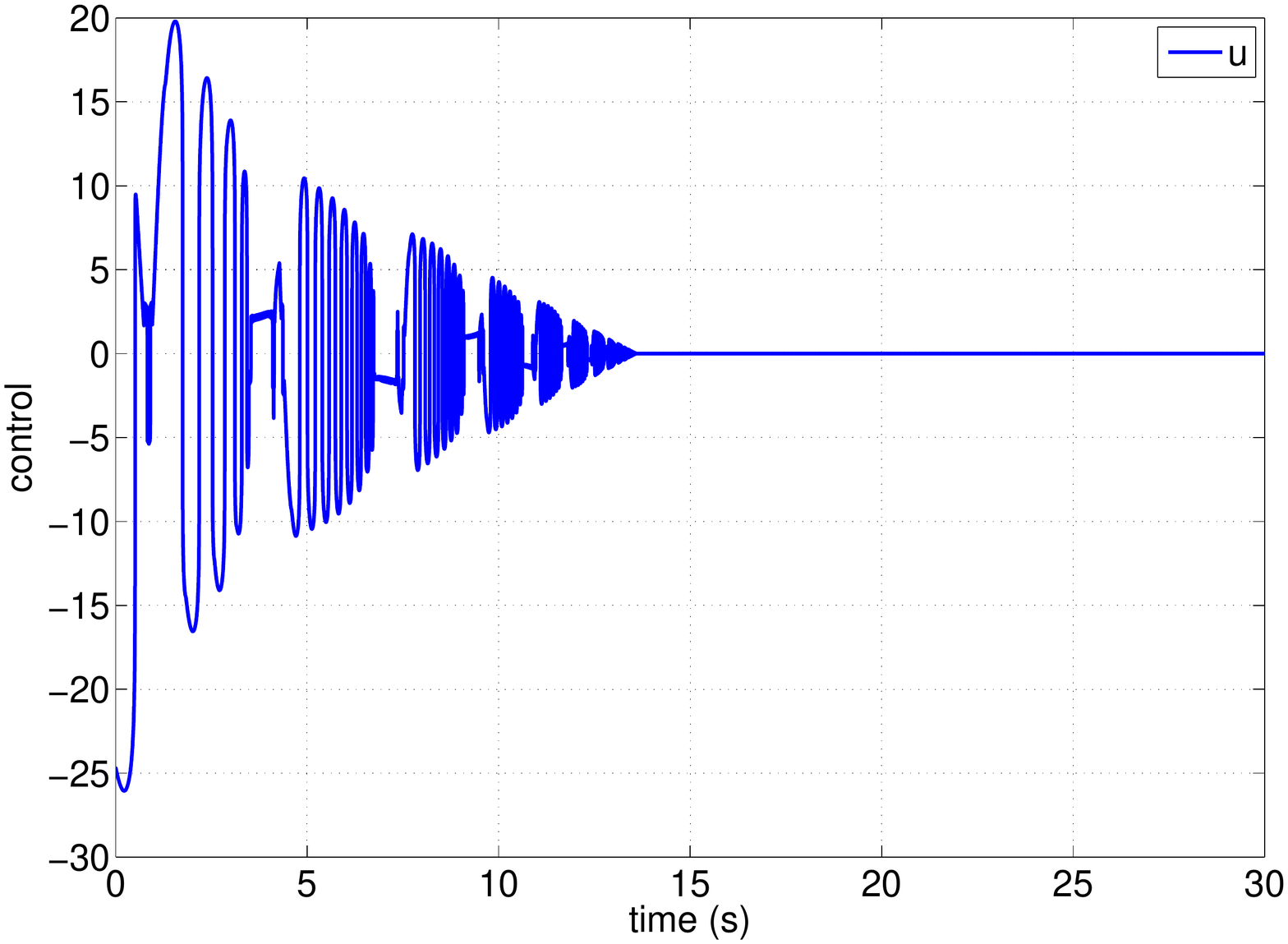}
    \label{u}
}\hspace{-1cm}
\subfigure[$z_1$, $z_2$, $z_3$ and $z_4$]{
    \includegraphics[width= 7 cm, height = 3.4 cm]{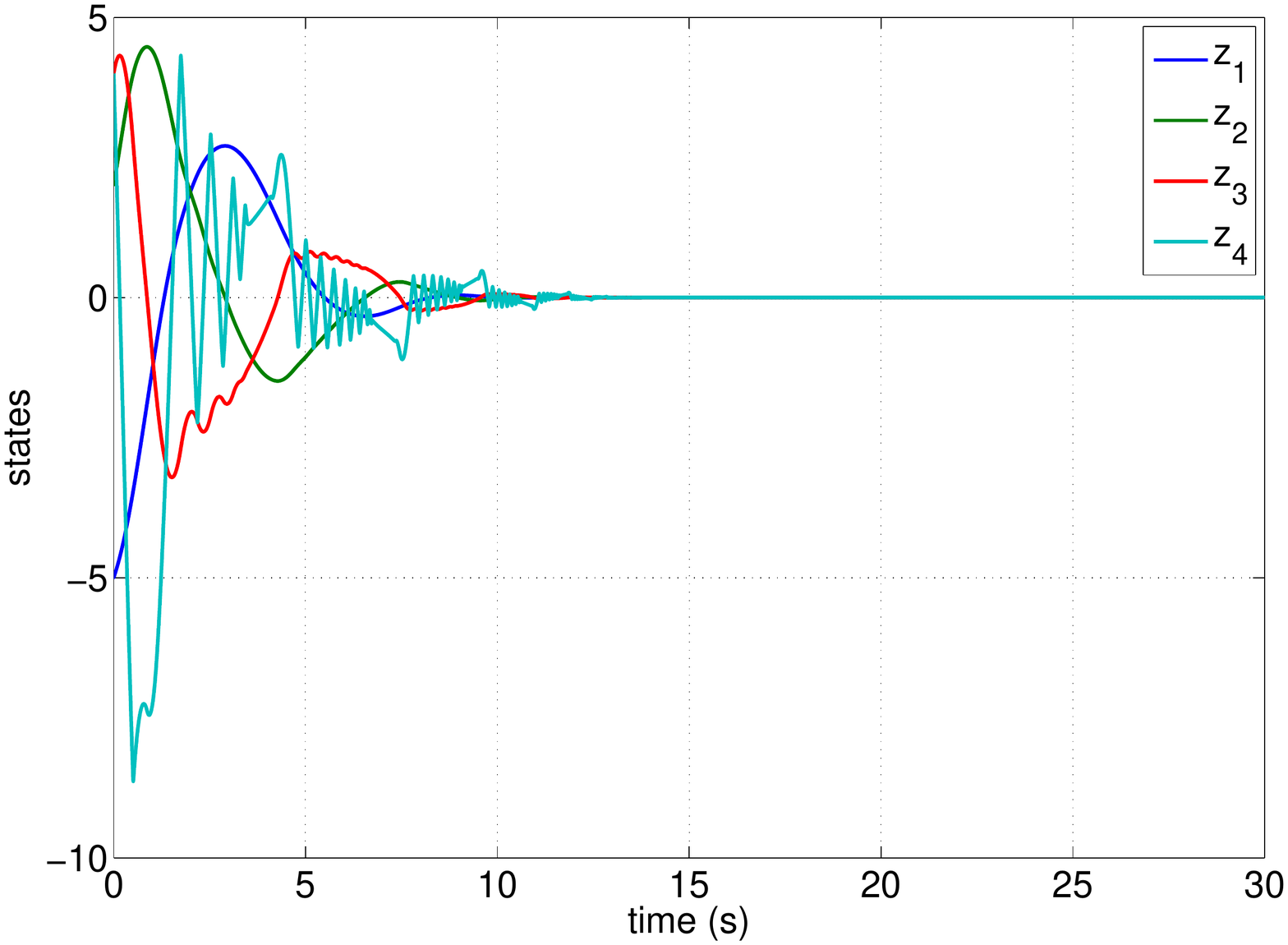}
    \label{s}
}%\hspace{-0.5cm}
%\subfigure[$\frac{\partial V_1}{\partial z_4}$ versus time ($s$).]{
    %\includegraphics[width= 8 cm, height = 3.6 cm]{D.pdf}
    %\label{d}
%}\hspace{-1cm}
\subfigure[Integral action $\int \partial_4 V_1 dt$]{
    \includegraphics[width= 7 cm, height = 3.4 cm]{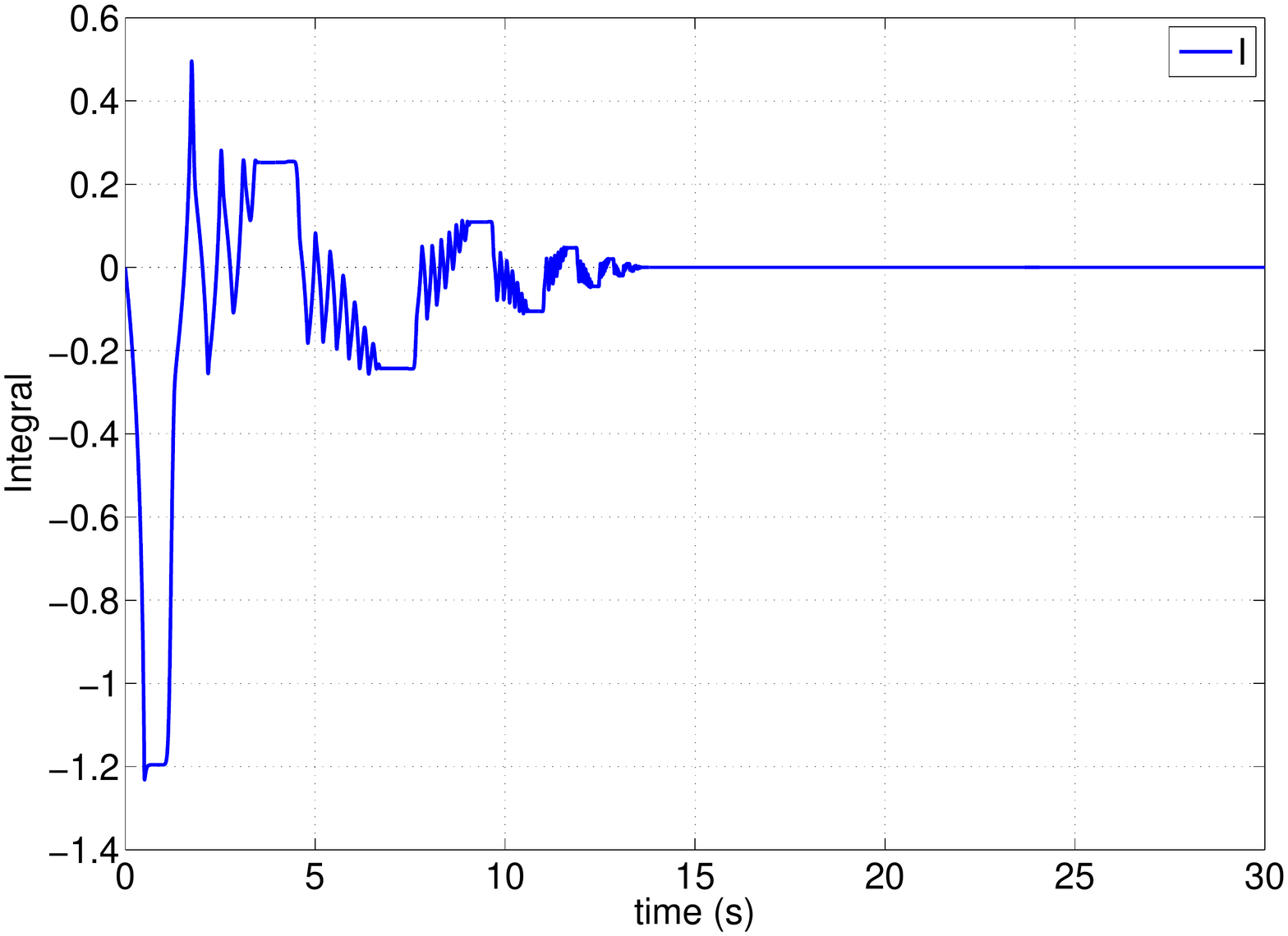}
    \label{i}
		}
		\caption{Pure integrator chain with integral action (with Hong's controller)}
\label{su0}
\end{figure}

\begin{figure}[htbp!]
\centering
\subfigure[Control law $u$]{
    \includegraphics[width= 7 cm, height = 3.4 cm]{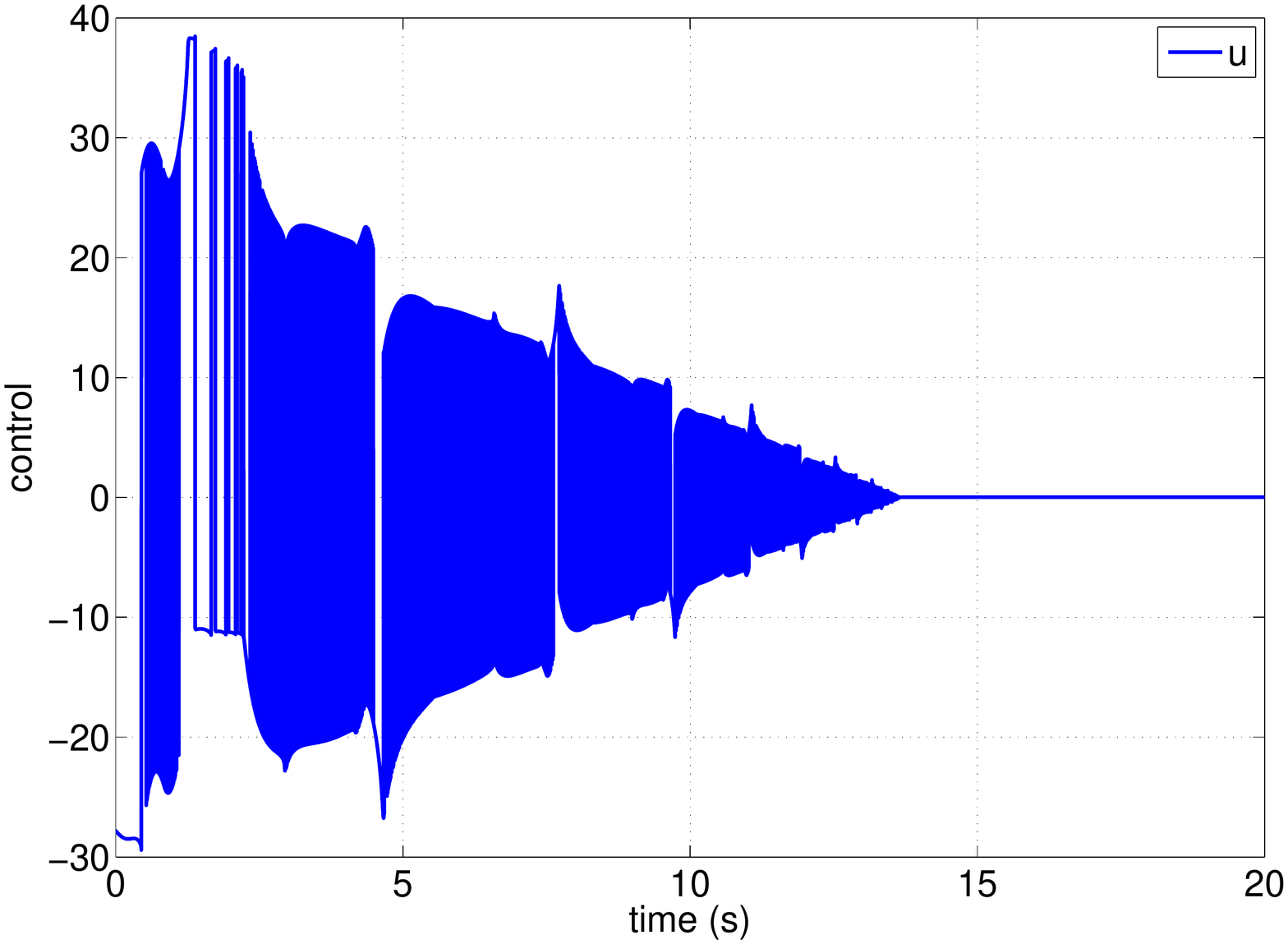}
    \label{u_qq}
}\hspace{-1cm}
\subfigure[$z_1$, $z_2$, $z_3$ and $z_4$]{
    \includegraphics[width= 7 cm, height = 3.4 cm]{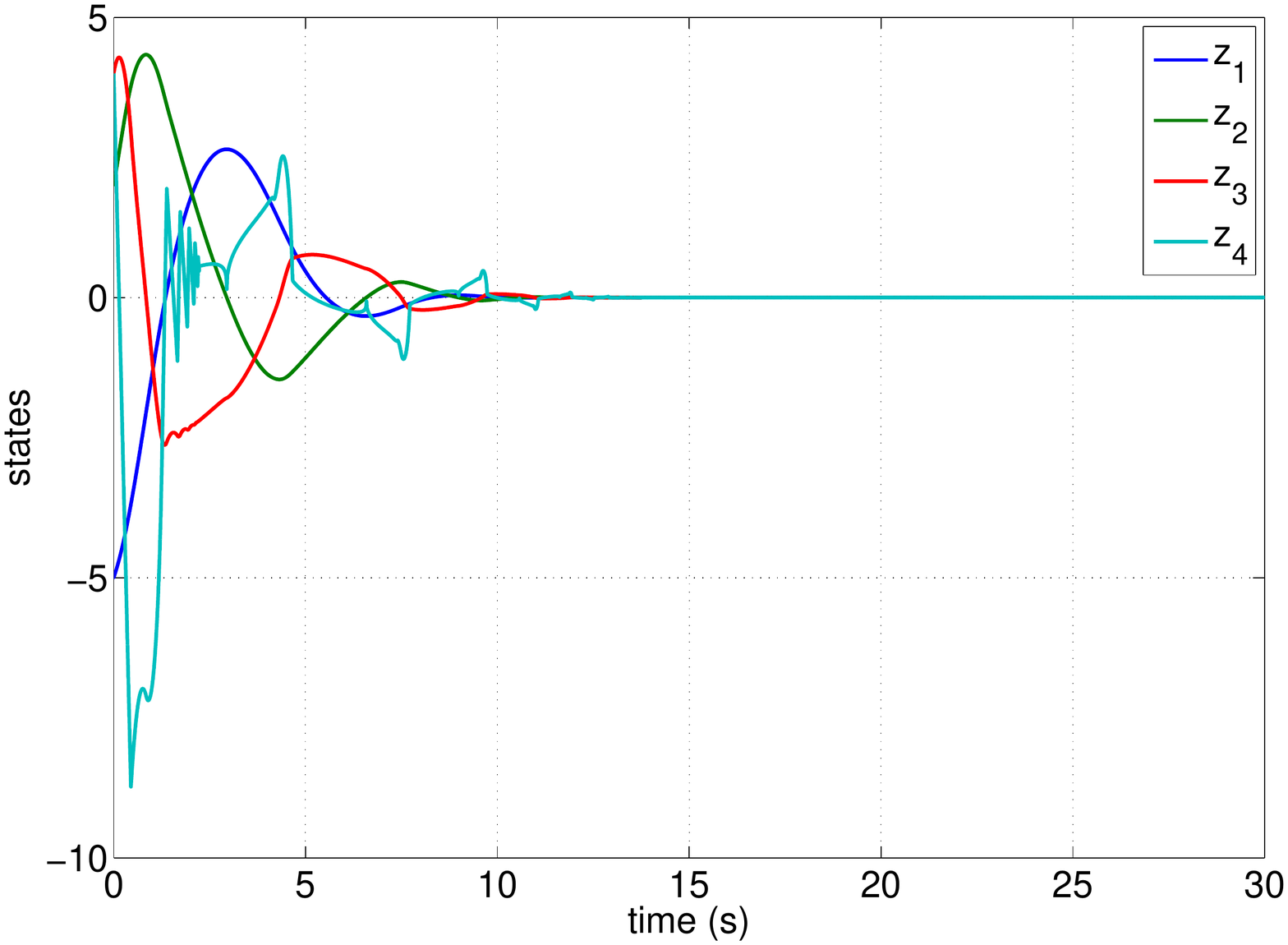}
    \label{s_qq}
}%\hspace{-0.5cm}
%\subfigure[$\frac{\partial V_1}{\partial z_4}$ versus time ($s$).]{
    %\includegraphics[width= 8 cm, height = 3.6 cm]{D_qq.pdf}
    %%\label{d_qq}
%}\hspace{-1cm}
\subfigure[ Integral action $I = \int \partial_4 V_1 dt$]{
    \includegraphics[width= 7 cm, height = 3.4 cm]{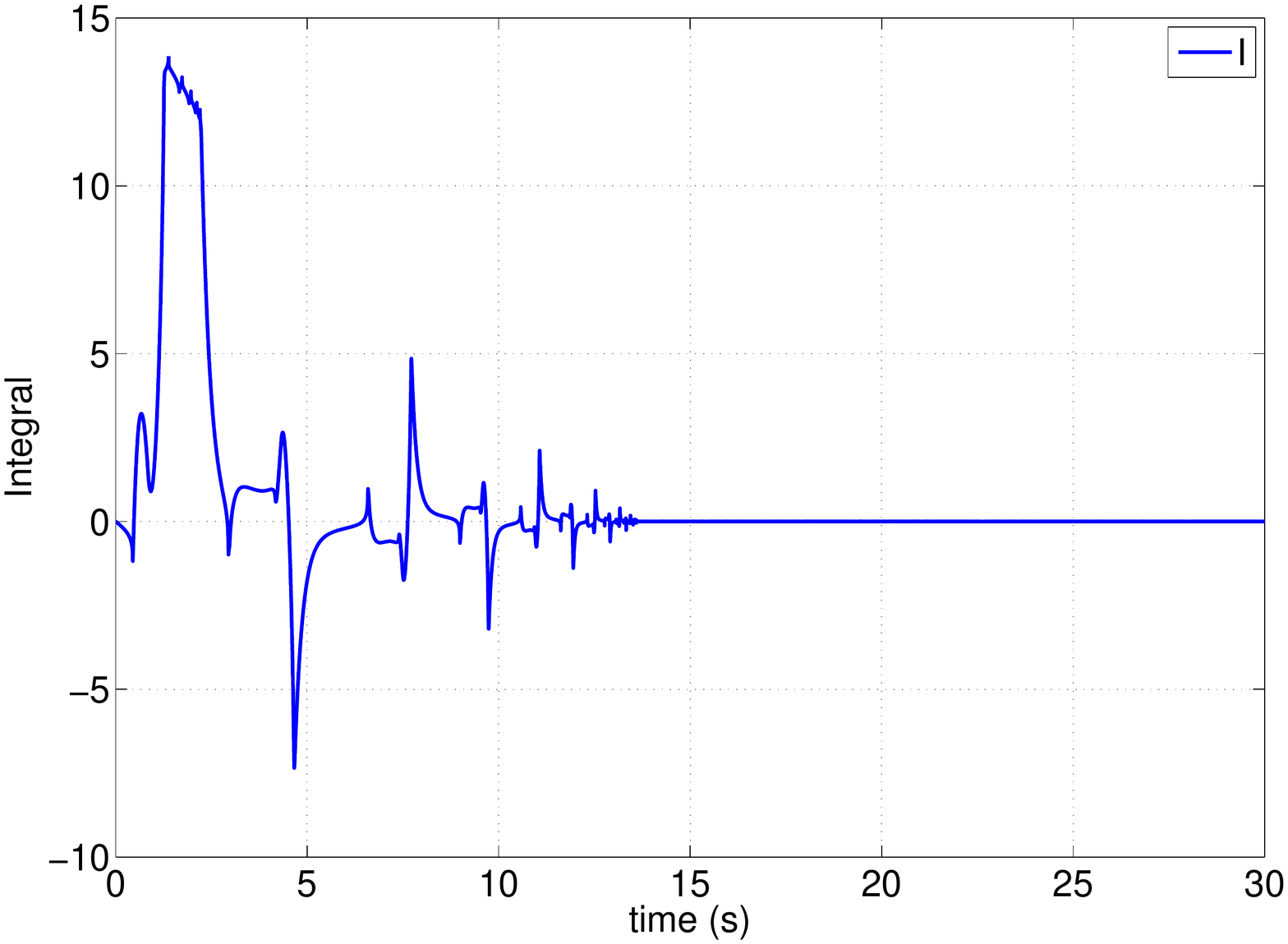}
    \label{i_qq}
		}
		\caption{Pure integrator chain with integral action (with modified Hong's controller)}
%\label{su0}
\end{figure}

In this subsection, we show the performance of HOST for the pure integrator chain for $u = k_P u_0 - k_I\int \partial_4V_1dt$.
The simulation parameters related to $u_0$ with the initial condition are tuned as in the previous subsection. The gains $k_P$ and $k_I$ are chosen as  $k_P=k_I = 1$.
The state convergence is presented in Figure \ref{s} and \ref{s_qq} for the Hong's controller
and the modified Hong's controller given in Figure \ref{u} andd Figure \ref{u_qq} respectively.
%In Figure \ref{d}, we remark that $\frac{\partial V_1}{\partial z_4}$ is continuous before the convergence to zero of the state.
Figures \ref{i} and \ref{i_qq} show the continuous integrator action wich vanishes to zero as there is no perturbation to compensate.

%simulation avec action integral
%%%%%%%%%%%%%%%%%%%%%%%%%%%%%%%%%%%%%%%%%%%%%%%%%%%%%%%%%%%%%%%%%%%%%%%%%%%%%%%%%
\subsection{Stabilisation of perturbed integrator chain by HOST - $\varphi\ne 0,\ \gamma \ne 1$}
We now consider the case of a perturbed system with perturbations $\varphi,\ \gamma$ defined as
\beq
		\varphi(t) &=& sin(t),\\
		\gamma(t) &=& 3 + \frac{1}{2} sin(0.5t),
\eeq
%$\varphi(t) = sin(t)$ and $|\dot\varphi| \le 1$. Note that $\varphi$ is bounded as well. $\gamma$ is given as $\gamma(t) = 3 + \frac{1}{2} sin(0.5t)$. 
Clearly, $\varphi$ is bounded and globally Lipschitz, as well as $\gamma$ which is in addition positive.

\begin{figure}[htbp!]
\centering
\subfigure[Control law $u$]{
    \includegraphics[width= 7 cm, height = 3.4 cm]{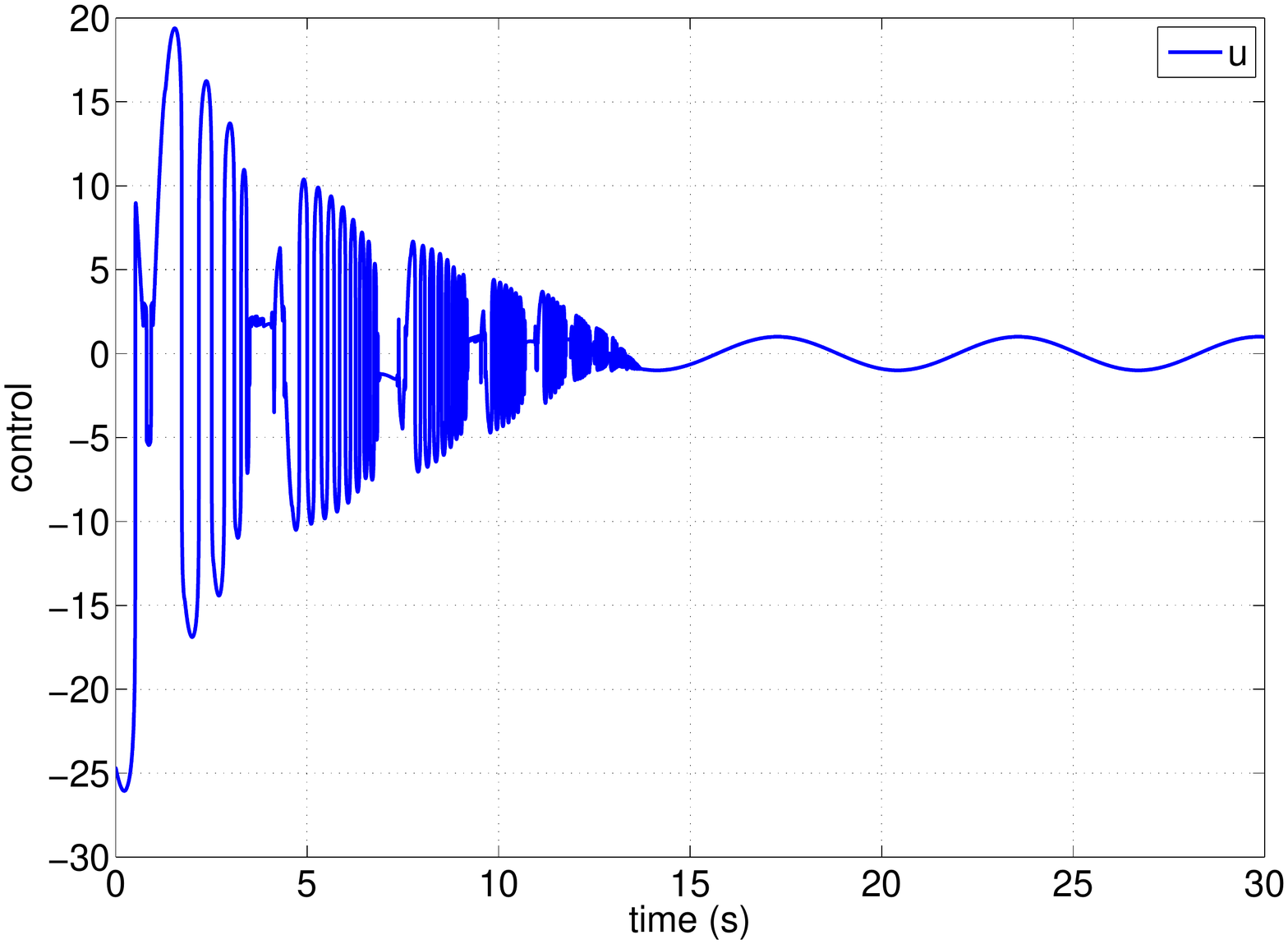}
    \label{up}
}\hspace{-1cm}
\subfigure[$z_1$, $z_2$, $z_3$ and $z_4$]{
    \includegraphics[width= 7 cm, height = 3.4 cm]{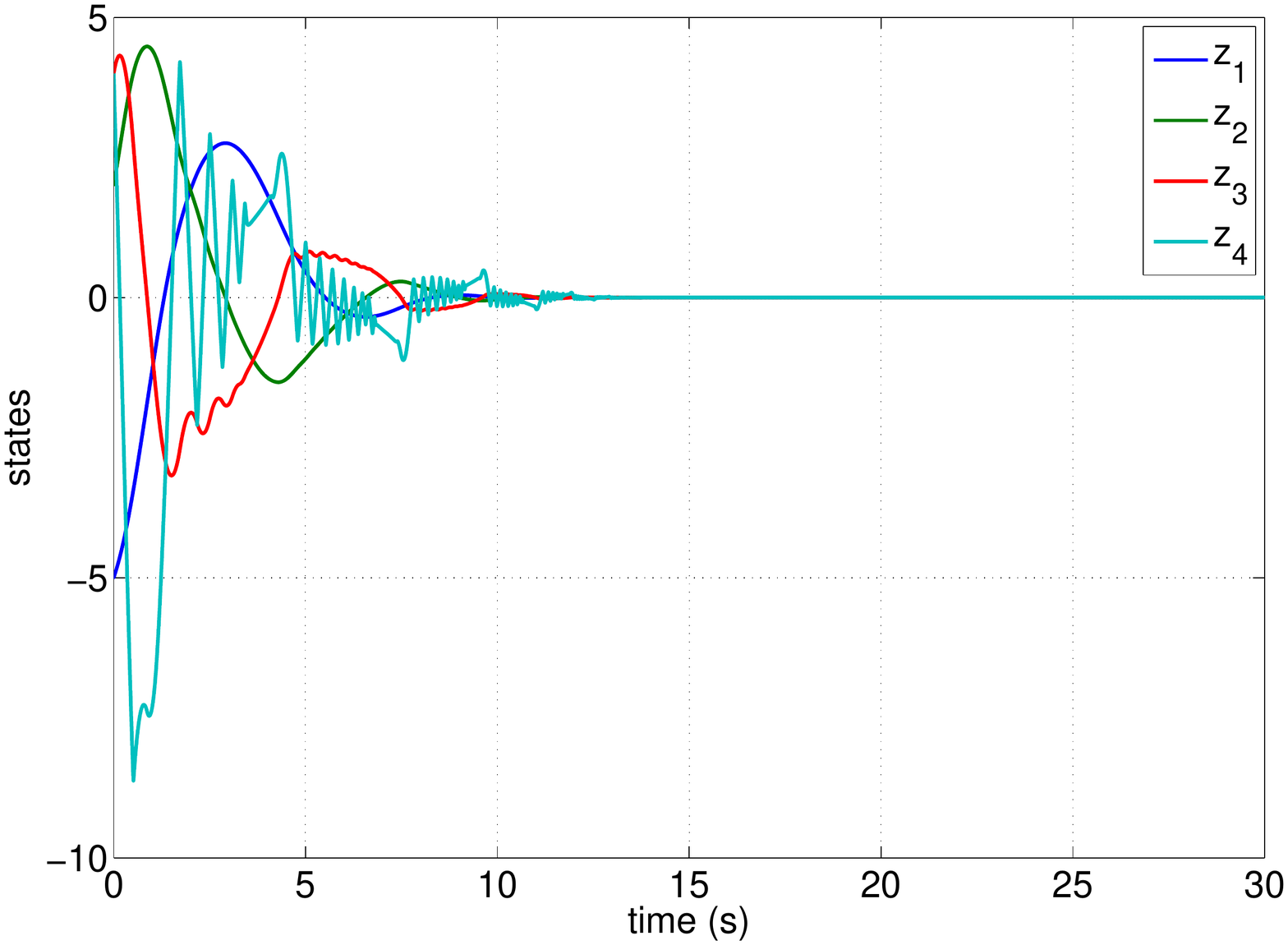}
    \label{sp}
%}\hspace{-0.5cm}
%\subfigure[$\frac{\partial V_1}{\partial z_3}$ versus time ($s$).]{
    %\includegraphics[width= 8 cm, height = 3.6 cm]{Dp.pdf}
    %\label{dp}
}\hspace{-1cm}
\subfigure[ Integral action $I = \int \partial_4 V_1 dt$]{
    \includegraphics[width= 10 cm, height = 3.4 cm]{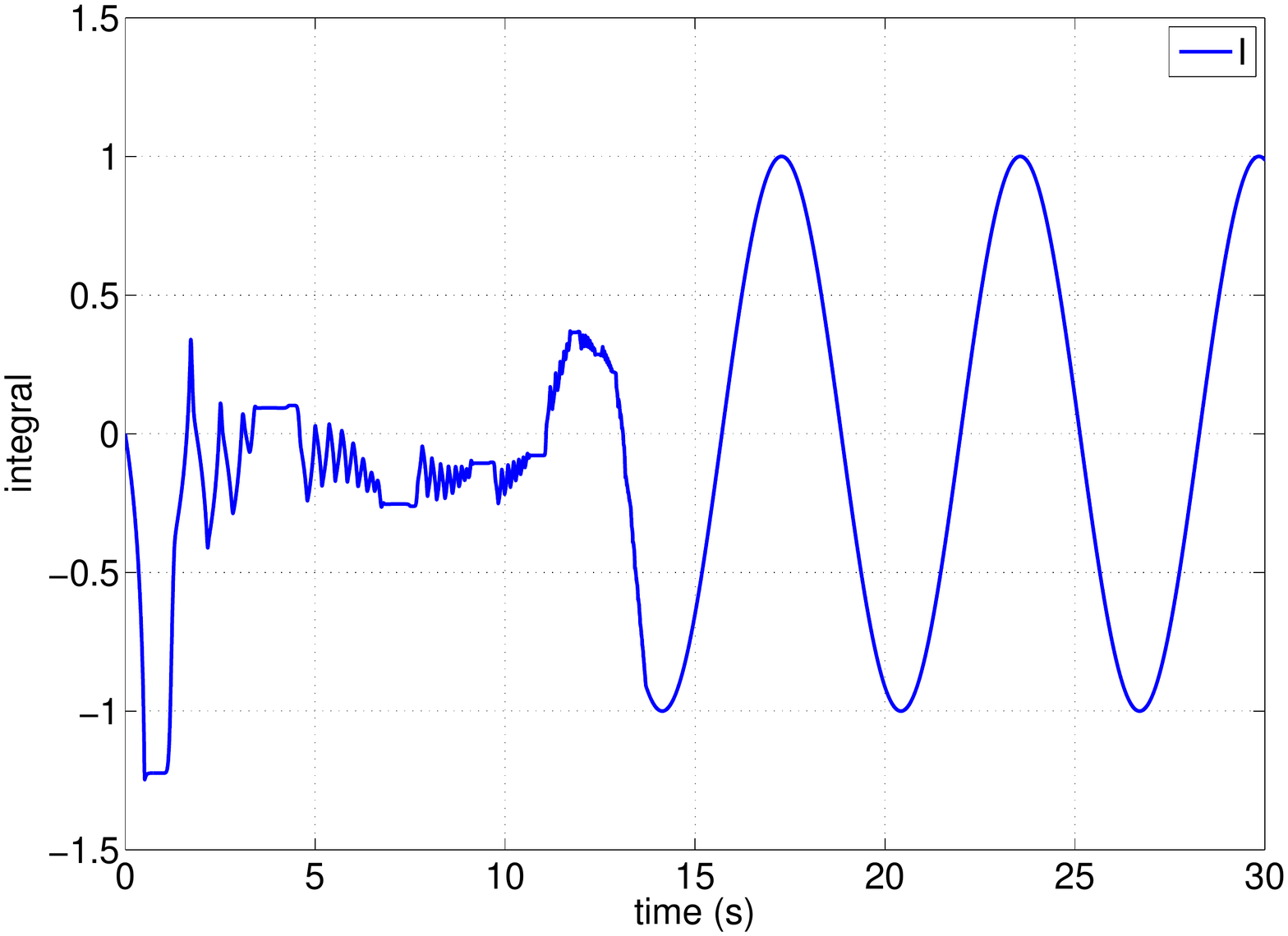}
    \label{ip}
}
\caption{Perturbed integrator chain with integral action (with Hong's controller)}
\label{s_u}
\end{figure}

\begin{figure}[ht]
\centering
\subfigure[Control law $u$]{
    \includegraphics[width= 7 cm, height = 3.4 cm]{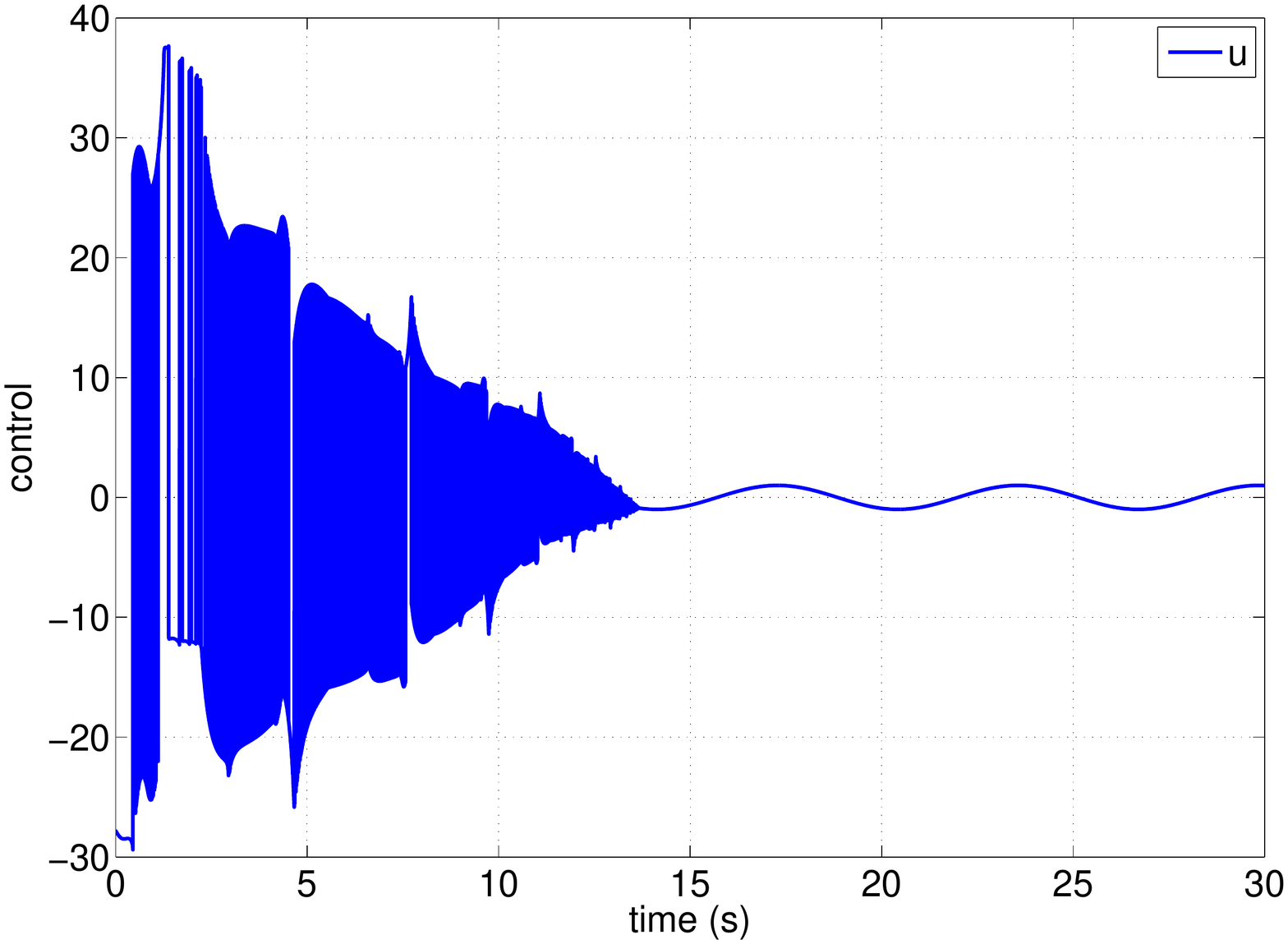}
    \label{up_qq}
}\hspace{-1cm}
\subfigure[$z_1$ and $z_2$ and $z_3$ and $z_4$]{
    \includegraphics[width= 7 cm, height = 3.4 cm]{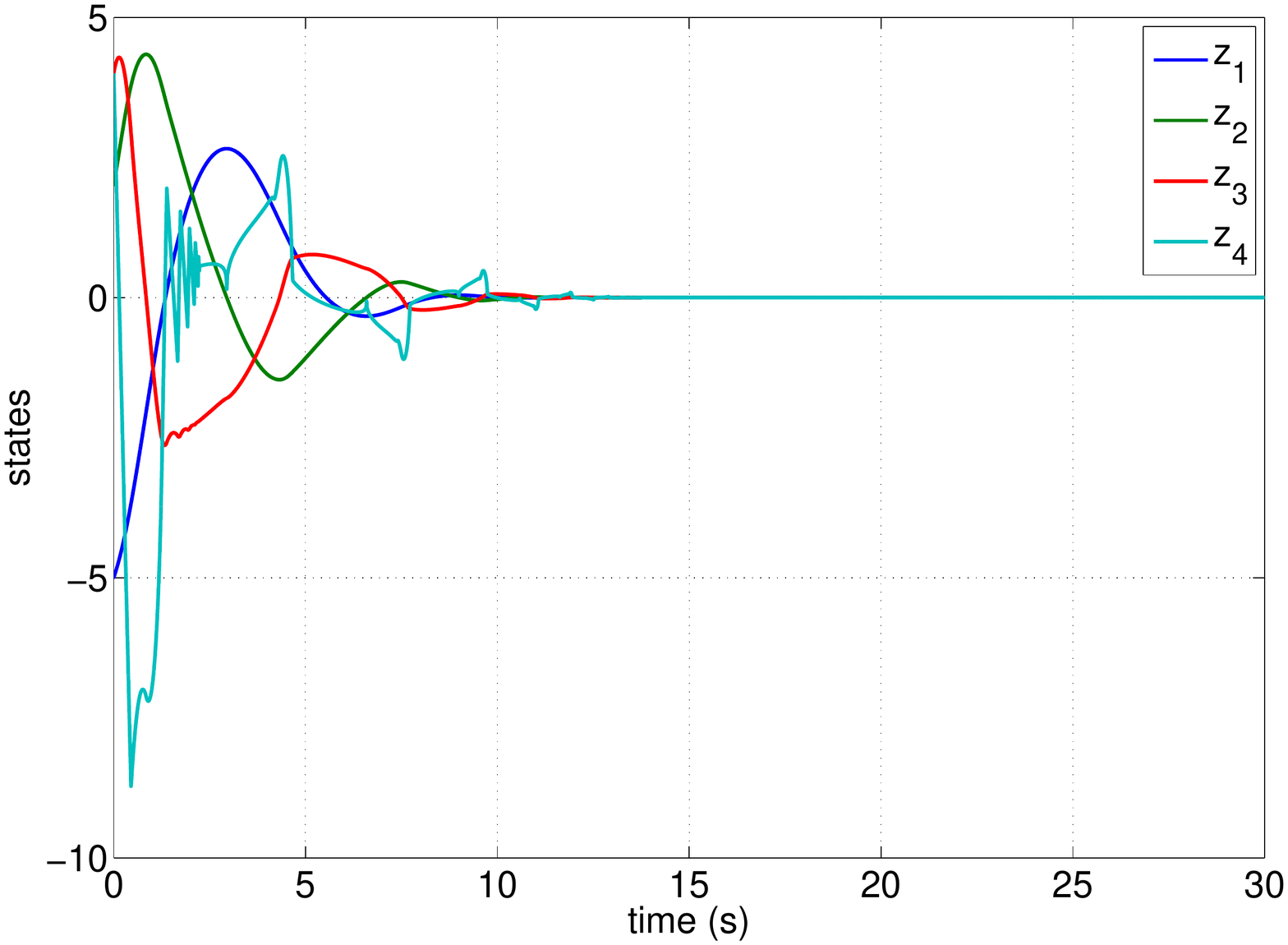}
    \label{sp_qq}
%}\hspace{-0.5cm}
%\subfigure[$\frac{\partial V_1}{\partial z_3}$ versus time ($s$).]{
    %\includegraphics[width= 8 cm, height = 3.6 cm]{Dp.pdf}
    %\label{dp}
}\hspace{-1cm}
\subfigure[ Integral action $ I = \int \partial_4 V_1 dt$]{
    \includegraphics[width= 10 cm, height = 3.4 cm]{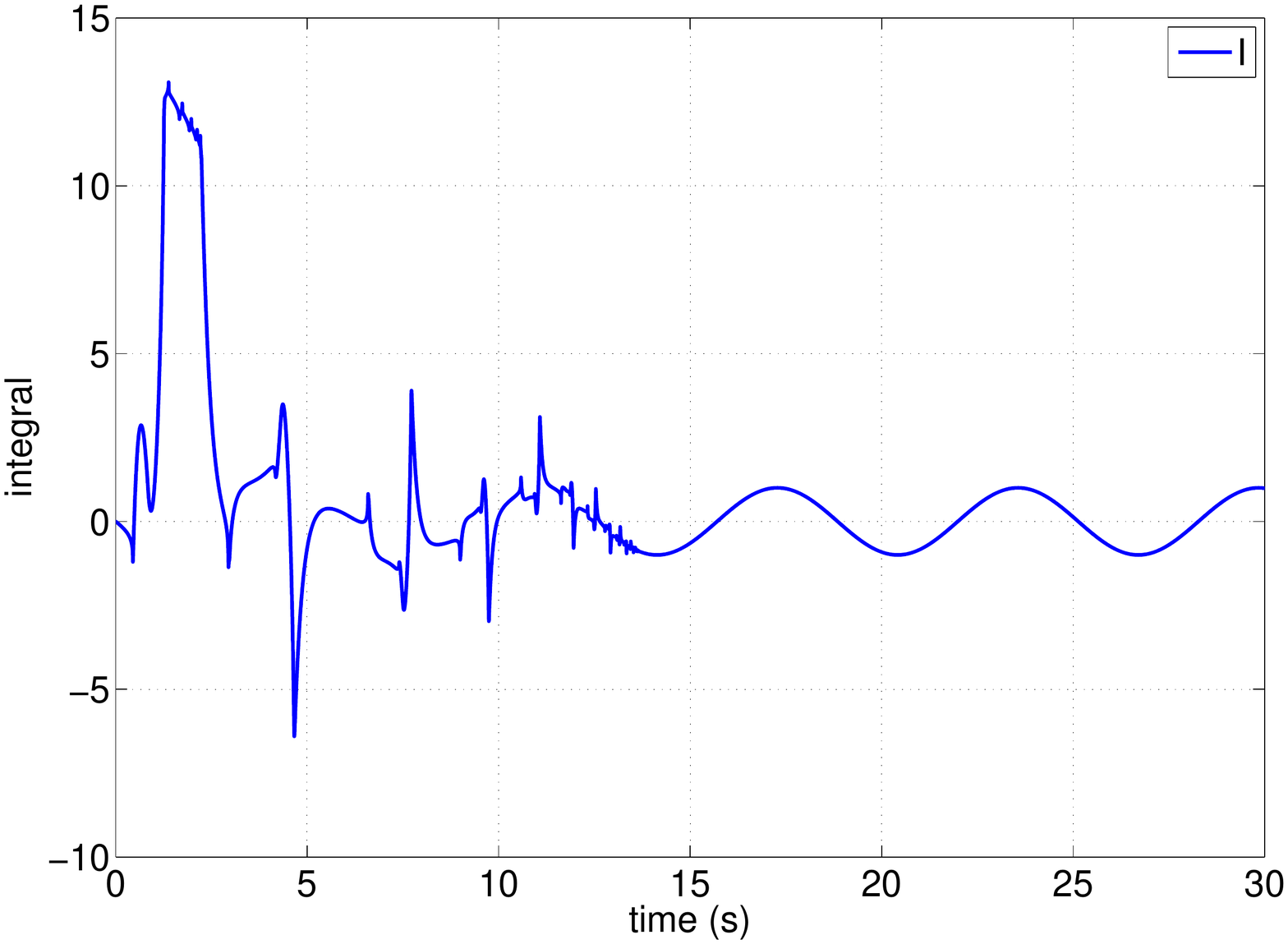}
    \label{ip_qq}
}
\caption{Perturbed integrator chain with integral action (Modified Hong's controller)}
\label{s_u_qq}
\end{figure}

The result is similar to the previous cases. However the controller acts in order to compensate the perturbation and we can see clearly in Figure \ref{up} and Figure \ref{ip} for the Hong's controller that $u(t) = -\frac{\varphi(t)}{\gamma(t)}$ after convergence to zero of the state. Similar results are obtained in the case of the modified Hong's controller in Figure \ref{up_qq} and Figure \ref{ip_qq}.

		%simulation du systeme perturbé

		%simulation du systeme perturbé

%%%%%%%%%%%%%%%%%%%%%%%%%%%%%%%%%%%%%%%%%%%%%%%%%%%%%%%%%%%%%%%%%%%%%%%%%%%%%%%%%%%%%%%%%%%%%%%%%%%%%%%%%%%%

%%%%%%%%%%%%%%%%%%%%%%%%%%%%%%%%%%%%%%%%%%%%%%%%%%%%%%%%%%%%%%%%%%%%%%%%%%%%%%%%%%%%%%%%%%%%%%%%%%%%%%%%%%%
%\pagebreak
\section{Conclusion}
In this paper we propose a general approach to design a continuous controller for a perturbed chain of integrators of arbitrary order generalizing the well-known supertwisting algorithm provided in \cite{Kamal2014} for integrator chain of length one and two. We have first designed a controller for the pure chain of integrators using a geometric condition inspired from \cite{Laghrouche_CST, Harmouche_CDC12} and we have proved convergence in finite time for the corresponding closed-loop system thanks to the explicit construction of strict homogeneous Lyapunov function for an extended system. As for the perturbed chain of integrators, we partially solve the complete problem by using homogeneity arguments applied to an extended differential inclusion. Future work consists of addressing the general case of a perturbed chain of integrators.% and of designing controllers only continous at the origin having better stabilization performances. %In order to choose efficiently the coefficients of the controllers in terms of the problem bounds and to get an estimate of the convergence time, it would be of great help to have a global strict Lyapunov function in the state variable and the supertwisting variable.

%%%%%%%%%%%%%%%%%%%%%%%%%%%%%%%%%%%%%%%%%%%%%%%%%%%%%%%%%%%%%%%%%%%%%%%%%%%%%%%%%%%%%%%%%%%%%%%%%%%%%%%%%%%
%\pagebreak
%\newpage
%\newpage
\bibliography{References_Stabilisation_2015_09_21}
\end{document}